\documentclass[12pt]{article}
\usepackage{color}
\usepackage{amsmath, amssymb, graphicx}
\begin{document}
\centerline{\bf Dynamic Multiplier Ideal Sheaves and}
\centerline{\bf the Construction of Rational Curves in Fano Manifolds}

\bigbreak\centerline{\it Dedicated to Professor Christer Kiselman} \bigbreak
\centerline{Yum-Tong Siu\ %
\footnote{Partially supported by a grant from the National Science
Foundation. Written for the Festschrift in honor of Professor Christer Kiselman
} }

\bigbreak

\bigbreak\noindent{\bf Introduction.}  Multiplier ideal sheaves were introduced by Kohn [Kohn1979] and Nadel [Nadel1990] to identify the location and the extent of the failure of crucial estimates.  Such multiplier ideal sheaves are defined by a family or a sequence of inequalities instead of a single inequality.  In Kohn's definition there is one inequality for every test function (or test form) and the multiplier has to make all the inequalities hold for all the test functions (or test forms) at the same time.  Nadel's definition is designed for the continuity method for the problem of the existence of K\"ahler-Einstein metrics on Fano manifolds.  Nadel's multiplier has to make the uniform finiteness of the integral from the crucial zero-order estimate hold for the entire sequence of perturbations of K\"ahler potentials occurring in the closed part of the continuity method.

\medbreak The notion of multiplier ideal sheaves used in algebraic geometry involves only one single inequality in the definition of a multiplier.  For a local plurisubharmonic function $\varphi$ on a domain $G$ in ${\mathbb C}^n$ the multiplier ideal sheaf ${\mathcal I}_\varphi$ used in algebraic geometry consists of all holomorphic function germs $f$ on $G$ such that $\left|f\right|^2e^{-\varphi}$ is locally integrable.  Only a single inequality $\int\left|f\right|^2e^{-\varphi}<\infty$ is used in characterizing $f\in{\mathcal I}_\varphi$.

\medbreak On the other hand, for Nadel's multiplier ideal sheaves a sequence of plurisubharmonic functions $\varphi_{t_\nu}$ for $\nu\in{\mathbb N}$ is used and his multiplier ideal sheaf consists of all holomorphic function germs $f$ satisfying
$\sup_{\nu\in{\mathbb N}}\int\left|f\right|^2e^{-\varphi_{t_\nu}}<\infty$, when he uses the multiplier ideal sheaf to handle the situation of closedness in the continuity method as $t_\nu\to t_*$ with $\nu\to\infty$.  Nadel's definition of a multiplier ideal sheaf uses a sequence of inequality with a uniform bound to characterize a multiplier in it.

\medbreak To emphasize the fundamental difference in these two definitions of multiplier ideal sheaves, we refer to the multiplier ideal sheaf used in algebraic geometry involving only one single inequality as a {\it static multiplier ideal sheaf} and refer to the multiplier ideal sheaf in the sense of Nadel involving a sequence of inequalities a {\it dynamic multiplier ideal sheaf}.  Of course, a static multiplier ideal sheaf in algebraic geometry is a special case of a dynamic multiplier ideal sheaf in the sense of Nadel when every term of the sequence $\varphi_{t_\nu}$ is equal to a fixed $\varphi$.

\medbreak A multiplier ideal sheaf in the sense of Kohn involves a family of inequalities parametrized by the collection of test functions (or test forms) and is also a dynamic multiplier ideal sheaf instead of a static multiplier ideal sheaf.

\medbreak The dynamic nature of dynamic multiplier ideal sheaves such as those in the sense of Nadel is specifically designed to terminate or stabilize a sequence of processes such as preventing a sequence of numbers or functions from increasing without bounds.  This powerful feature is no longer found in static multiplier ideal sheaves used in algebraic geometry.

\medbreak In the analytic proof of the finite generation of the canonical ring for a compact complex manifold of general type, dynamic multiplier ideal sheaves are used  (together with the notion of deviation from sufficient ampleness instead of minimum centers of log canonical singularities as deviation from freeness) [Siu2006, Siu2007, Siu2008, Siu2008a].  That is the reason why the infinite process of blowing up the base-point set to hypersurfaces of normal crossing can be terminated in the analytic proof.  The dynamic nature of dynamic multiplier ideal sheaves enables us to terminate such an infinite process.
Actually such ingredients of dynamic multiplier ideal sheaves and deviation from sufficient ampleness are already used without explicit mention in the technique of pluricanonical extensions introduced for the confirmation of the conjecture on the deformational invariance of plurigenera.  The technique of pluricanonical extension from analysis is also a crucial step in the algebraic geometric approach to the problem of the finite generation of the canonical ring for a compact complex manifold of general type [Birkan-Cascini-Hacon-McKernan2006].  The advantage of the analytic proof of the finite generation of the canonical ring is that its use of dynamic multiplier ideal sheaves explains completely transparently why the infinite blow-up process terminates and why the argument works (see [Siu2008, Remark(1.3.1)] and [Siu2008a, (1.2)(H)]).

\medbreak Actually the use of the semi-continuity of multiplier ideal sheaves to treat the freeness of the Fujita conjecture in [Angehrn-Siu1995] also stemmed from dynamic multiplier ideal sheaves though there was no explicit mention of it there.

\medbreak In this note we will discuss and explain the historic evolution of the notion of multiplier ideal sheaves, especially the interpretation from the viewpoint of destabilizing subsheaves in the context of terminating or bounding an infinite process.  We will start out with Kohn's subelliptic multipliers and explain its relation with Nadel's multipliers.  We will use the construction of Hermitian-Einstein metrics for stable vector bundles to heuristically illustrate the viewpoint of interpreting multiplier ideal sheaves as destabilizing subsheaves.

\medbreak We will also discuss the approach of constructing rational curves in Fano manifolds by using dynamic multiplier ideal sheaves and {\it singularity-magnifying} complex Monge-Amp\`ere equations.  This approach is still under development with details in the process of being worked out.  We will indicate where details still need to be worked out.  This part of our note is only a presentation of our approach together with various techniques and ideas which we have been developing for it.  A complete analytic proof of the existence of rational curves in Fano manifolds is not yet available.  Our approach is presented here to open up a new direction and to introduce a new area of research in the interface between several complex variables and algebraic geometry.

\medbreak The only known method of constructing rational curves in Fano manifolds is the bend-and-break method of Mori [Mori1979] using the method of characteristic $p>0$.  For three decades it has been a challenge to complex geometers and global analysts to find a way to prove the existence of rational curves in Fano manifolds by using analytic methods without involving characteristic $p>0$.  The only result relevant for this problem obtained by methods of complex geometry is the use of energy-minimizing harmonic maps in [Siu-Yau1980] to produce rational curves in compact complex manifolds of positive bisectional curvature, but such a technique is useless for the problem of constructing rational curves in Fano manifolds.  The approach presented here of using dynamic multiplier ideal sheaves and singularity-magnifying complex Monge-Amp\`ere equations is an endeavor to remove this three-decade-old thorn on the side of analysts.

\medbreak When one uses the theorem of Hirzebruch-Riemann-Roch [Hirzebruch1966] and multi-valued holomorphic anticanonical sections to obtain static multiplier ideal sheaves to construct rational curves in Fano manifolds, one encounters the problem of insufficient size of the relevant Chern classes, just like in the Mori's bend-and-break technique of deforming a curve with two points fixed before his introduction of the method of characteristic $p>0$.  Mori's use of the method of characteristic $p>0$ enables him to increase the relevant Chern classes to the necessary size.  Our approach of using dynamic multiplier ideal sheaves and singularity-magnifying complex Monge-Amp\`ere equations to produce destabilizing subsheaves serves the same function of removing the limitation imposed by the insufficient size of the relevant Chern classes. This in some way corresponds to the r\^ole of Mori's use of the method of characteristic $p>0$.  For the benefit of analysts reading this note in an appendix we will highlight the key points of Mori's argument for comparison with our approach. At a very loose philosophical level the use of the Monge-Amp\`ere equation singularity-magnifying complex Monge-Amp\`ere equations to produce destabilizing subsheaves is dual to Mori's method of deforming complex curves.  A subvariety in a complex manifold can be defined by a map from a compact complex space to the manifold or by a coherent ideal sheaf or more generally a coherent subsheaf.  Mori's approach uses the deformation of a holomorphic map and in our approach we use the deformation of a coherent subsheaf giving rise to a destabilizing subsheaf.

\medbreak We would like to inject here another remark about the difference between producing a coherent subsheaf by using Chern classes and the theorem of Hirzebruch-Riemann-Roch [Hirzebruch1966] and producing a destabilizing subsheaf.  This kind of difference was already used in the construction of holomorphic sections of ample vector bundles, especially the construction of holomorphic jet differentials for hyperbolicity problems, first by Miyaoka [Miyaoka1983], and then by Schneider-Tancredi [Schneider-Tancredi1988], and by Lu-Yau [Lu-Yau1990].

\medbreak Besides the use of Nadel's vanishing theorem the relation between destabilizing subsheaves and rational curves can be seen also from the phenomenon that though the tangent bundle of ${\mathbb P}_n$ for $n\geq 2$ is stable its restriction to a minimal rational curve of ${\mathbb P}_n$ is not.  The restriction of a stable vector bundle to a curve of appropriate genericity is stable but minimal rational curves do not belong to such a class of curves so far as the tangent bundle of ${\mathbb P}_n$ is concerned (see (2.7) below).

\medbreak Various parts of the content of this note has been presented in a number of recent conferences.  There have been many requests from the participants of the conferences for computer files used in the presentations.  One of the reasons for making this note available is to respond to such requests.

\bigbreak\noindent{\bf \S1. Historic Evolution of Multiplier Ideal Sheaves.}

\bigbreak\noindent(1.1) {\it Kohn's Subelliptic Multipliers
for the Complex Neumann Problem.}  The setting is a bounded domain $\Omega$ in ${\bf C}^n$ with smooth
weakly pseudoconvex boundary defined by $r<0$ with $dr$ being nowhere
zero on the boundary $\partial\Omega$ of $\Omega$. Here weakly pseudoconvex boundary means that
$\sqrt{-1}\,\partial\bar\partial r|_{T^{(1,0)}_{\partial\Omega}}\geq 0$.  The problem is to study the following regularity question: given a smooth $(0,1)$-form $f$ on $\bar\Omega$ with $\bar\partial f=0$, whether the solution of
$\bar\partial u=f$ on $\Omega$ with $u$ perpendicular to all holomorphic functions on
$\Omega$ is smooth on $\bar\Omega$.

\medbreak A sufficient condition for regularity is the following subelliptic estimate at every boundary point.  For $P\in\partial \Omega$ there exist some
open neighborhood $U$ of $P$ in ${\bf C}^n$ and positive numbers
$\epsilon$ and $C$ satisfying
$$
\||g|\|_\epsilon^2\leq C\left(\|\bar\partial g\|^2+\|\bar\partial^*
g\|^2+\|g\|^2\right)
$$
for every $(0,1)$-form $g$ supported on $U\cap\bar\Omega$ which is in the
domain of $\bar\partial$ and $\bar\partial^*$.  Here $\||\cdot|\|_\epsilon$ is the $L^2$ norm on $\Omega$
involving derivatives up to order $\epsilon$ in the boundary
tangential directions of $\Omega$, $\|\cdot\|$ is the usual $L^2$ norm on $\Omega$ without involving any
derivatives, and $\bar\partial^*$ is the actual adjoint of
$\bar\partial$ with respect to $\|\cdot\|$.

\medbreak The reason why some positive $\varepsilon$ is needed is that in applying a differential operator $D$ to both sides of $\bar\partial u=f$ to get estimates of the Sobolev norm of $u$ up to a certain order of derivatives in terms of that of $f$, an error term from the commutator of the differential operator $D$ and $\bar\partial$ occurs, which needs to be absorbed and one way to do the absorption is to use an estimate involving a Sobolev norm with derivative higher by some positive number $\varepsilon$.  This stronger Sobolev norm is used also to absorb the error term from partitions of unity or cut-off functions.

\medbreak The reason why only the tangential Sobolev norm $\||\cdot|\|_\epsilon$ is used is that we need to preserve the condition that $(0,1)$-form $g$ belongs to the domain of $\bar\partial^*$ (which means the vanishing of the complex-normal component at boundary points) by using only differentiation along the boundary tangential directions.  The missing estimate in the real-normal direction can be obtained from the complex-normal component of the equation $\bar\partial u=f$.

\medbreak The theory of multiplier ideal sheaves introduces multipliers into the most crucial estimate, which in this case is the subelliptic estimate.  Later in (1.4) Nadel's multiplier ideal sheaves will also be in like manner defined from the most crucial estimate in Nadel's setting.  For Kohn's setting here, a {\it subelliptic scalar multiplier} $F$ is a smooth function germ of ${\mathbb C}^n$ at $P$ such that the following subellitpic
estimate of some positive order $\varepsilon_F$ holds for any test $(0,1)$-form $g$
after replacing it by its product with $F$.
$$
\||Fg|\|_{\epsilon_{{}_F}}^2\leq C_{{}_F}\left(\|\bar\partial
g\|^2+\|\bar\partial^* g\|^2+\|g\|^2\right)
$$
for every test $(0,1)$-form $g$ described above.  The {\it multiplier ideal} $I_P$ at the boundary point $P$ is the ideal of all such subelliptic scalar multipliers $F$.

\medbreak A {\it subelliptic vector-multiplier} $\theta$ is a smooth $(1,0)$-form germ on ${\mathbb C}^n$ at $P$ such that the following
subellitpic estimate of some positive order $\varepsilon_\theta$ holds for any test
$(0,1)$-form $g$ after replacing it by its inner product $\theta\cdot g$ with $\theta$.
$$
\||\theta\cdot g|\|_{\epsilon_{{}_\theta}}^2\leq
C_{{}_\theta}\left(\|\bar\partial g\|^2+\|\bar\partial^*
g\|^2+\|g\|^2\right)
$$
for every test $(0,1)$-form $g$ described above.  The {\it multiplier module} $A_P$ at the boundary point $P$ is the module of all such subelliptic vector-multipliers $\theta$.

\bigbreak The most important part of the theory of Kohn's multiplier ideal sheaves is the following Kohn's Algorithm.

\begin{itemize}\item[(A)] {\it Initial Membership}.

\begin{itemize}\item[(i)] $r\in I_P$.

\item[(ii)] $\partial\bar\partial_j r$ belongs to
$A_P$ for every $1\leq j\leq n-1$ if $\partial r=d z_n$ at
$P$ for some local coordinate system $\left(z_1,\cdots,z_n\right)$,
where $\partial_j$ means $\frac{\partial}{\partial z_j}$.
\end{itemize}
\item[(B)] {\it Generation of New Members}.

\begin{itemize}\item[(i)]  If $f\in I_P$, then $\partial f\in A_P$.

\item[(ii)] If $\theta_1,\cdots,\theta_{n-1}\in
A_P$, then the coefficient of
$
\theta_1\wedge\cdots\wedge\theta_{n-1}\wedge\partial r$ is in
$I_P$. \end{itemize}\item[(C)] {\it Real Radical Property}.

\smallbreak\noindent If $g\in I_P$ and
$\left|f\right|^m\leq\left|g\right|$ for some positive integer $m$, then $f\in I_P$.
\end{itemize}
\medbreak Kohn's algorithm allows certain differential operators to
lower the vanishing order of multiplier ideals.  There are the following two limitations on using differentiation to reduce vanishing orders.  The first one is that only $(1,0)$-differentiation is
allowed.  The second one is that only determinants of coefficients of $(1,0)$-differentials (in the complex tangent space of the boundary)
from Cramer's rule can be used. Moreover, root-taking can be used to reduce vanishing orders.

\medbreak The goal of Kohn's algorithm is to produce the constant function $1$ as a subelliptic scalar multiplier under some appropriate geometric assumption on the boundary.  The geometric assumption is the following finite type condition formulated by D'Angelo [D'Angelo1979] and the goal is to verify Kohn's conjecture which is given below [Kohn1979].

\medbreak\noindent The {\it type} $m$ at a point $P$ of the
boundary of weakly pseudoconvex $\Omega$ is the supremum of the
normalized touching order
$$\frac{{\rm ord}_0\left(r\circ\varphi\right)}{{\rm ord}_0\varphi}$$
to $\partial\Omega$, of all local holomorphic curves
$\varphi:\Delta\to{\mathbb C}^n$ with $\varphi(0)=P$, where $\Delta$
is the open unit $1$-disk and ${\rm ord}_0$ is the vanishing order
at the origin $0$.  The domain $\Omega$ is of finite type if the supremum of the type of
every one of its boundary points is finite.

\medbreak\noindent{\it Kohn's Conjecture:} Kohn's algorithm
terminates for smooth weakly pseudoconvex domains of finite type
(with effectiveness involving type and order of subellipticity).

\medbreak\noindent Kohn's conjecture was solved for the real-analytic case without
effectiveness [Diederich-Fornaess1978].  A more geometric proof
of Proposition 3 in [Diederich-Fornaess1978] (which is the key step) is given in [Siu2007] where the geometric viewpoint delineates more the r\^ole played by the real-analytic assumption and the hurdle standing between generalizing the ineffective real-analytic case to the ineffective smooth case.  The effective termination of Kohn's algorithm is given in [Siu2007] for the case where $\Omega$ is a special domain in the sense of [Kohn1979] and how the techniques given there are to be extended to give a proof of the full Kohn conjecture with effectiveness is also described in it.

\medbreak Kohn's algorithm can be geometrically interpreted in terms of the usual Frobenius theorem on the integrability of a distribution of the linear subspace of the tangent space, which states that for an open subset $U\subset{\mathbb R}^m$ and a distribution of $k$-dimensional subspace $x\mapsto V_x\subset T_{{\mathbb R}^m}={\mathbb
R}^m$ of the tangent space $T_{{\mathbb
R}^m}$ of ${\mathbb R}^m$, the distribution $V_x$ is integrable ({\it i.e.} $V_x$ is the tangent
space of a family of $k$-folds in $U$) if and only if the distribution is closed under Lie bracket in the sense that $\left[V_x,V_x\right]\subset
V_x$ for all $x\in U$ or alternatively $d\omega_j=\sum_{\ell=1}^{m-k}\omega_\ell\wedge\eta_{j,\ell}$, where $\omega_1,\cdots,\omega_{m-k}$ are $1$-forms defining
$V_x$ and $\eta_{j,1},\cdots,\eta_{j,m-k}$ are some other $1$-forms.

\medbreak There are other weaker forms of integrability than the full integrability in Frobenius's theorem.  For example, in his 1909 paper on thermodynamics Carath\'eodory [Carath\'eodory1909] introduced the notion of the weaker notion of integrability along curves in the case of codimension $1$ with $k=m-1$.  He considered smooth curves $C$ whose tangents at the point $x$ belong to $V_x$.  Chow in 1939 generalized Carath\'eodory's weaker notion of integrability along curves to the case of a general $1\leq k\leq m-1$ [Chow1939].

\medbreak In interpreting Kohn's algorithm in terms of Frobenius's integrability theorem, we consider a notion of integrability even weaker than that of Carath\'eodory and Chow along curves.  We consider integrability over an Artinian subscheme.
For example, the ringed space
$\left(0,\,{\mathcal O}_{{\mathbb C}^n}\left/{\mathcal
I}\right.\right)$ with $\left({\mathfrak m}_{{\mathbb
C}^n,0}\right)^N\subset{\mathcal I}$ for some integer $N\geq 1$ is an Artinian
subscheme.  If the distribution $V_x$ is in an open subset $U$ of ${\mathbb R}^m\subset{\mathbb C}^n$ with $m\leq 2n$, the integrability over the Artinian subscheme $\left(0,\,{\mathcal O}_{{\mathbb C}^n}\left/{\mathcal
I}\right.\right)$ is the same as some corresponding jet of an
complex curve at $0$ is tangential to the distribution of tangent subspaces.

\medbreak For the interpretation of Kohn's algorithm we denote by $M$ the real hypersurface $M$ which is the boundary $r=0$ of $\Omega$ in ${\mathbb C}^n$ and consider the distribution $T^{\mathbb R}_M\cap JT^{\mathbb R}_M$ on $M$, where $J$ is the almost complex structure of $M$.  The usual full Frobenius integrability means that $M$ is Leviflat.
Integrability over an Artinian subscheme of high order means
some local holomorphic curve touching $M$ to high order at one point.  Finite type in the sense of D'Angelo means a limit on the order of the Artinian subscheme of
integrability.  Kohn's Algorithm is simply the condition, expressed in terms of differential forms defining the distribution, for limiting the order of an Artinian subscheme of integrability.  It is similar to the condition $d\omega_j=\sum_{\ell=1}^{m-k}\omega_\ell\wedge\eta_{j,\ell}$ for the usual Frobenius theorem, but points to the opposite direction.  The condition $d\omega_j=\sum_{\ell=1}^{m-k}\omega_\ell\wedge\eta_{j,\ell}$ means that when we differentiate $\omega_j$ to get $d\omega_j$, we do not get anything new, because the result $d\omega_j$ is already generated by $\omega_\ell$ for $1\leq\ell\leq m-k$. In contrast, Kohn's algorithm starts out with $\partial r$ which defines the distribution $T^{\mathbb R}_M\cap JT^{\mathbb R}_M$ and, when we take its differential $d\partial r=-\partial\bar\partial r$, we get something new and when we use Cramer's rule and other procedures, we keep on getting something new until we end up with the constant function $1$ as a subelliptic scalar multiplier.  The effectiveness involved in the procedure of getting finally the constant function $1$ places a limit on the the order of an Artinian subscheme of integrability.

\bigbreak\noindent(1.2) {\it Interpretation of Multiplier Ideal Sheaves in Terms of Rescaling and as Destabilizing Subsheaves.} Multipliers are introduced into crucial {\it a priori} estimates occurring in the solution or the regularity problem for a differential equation $Lu=f$.  For the regularity problem like the situation of the complex Neumann problem considered by Kohn, the {\it a priori} estimates will in general involve some stronger norms such as $\left|\left\|\cdot\right\|\right|_\varepsilon$.

\medbreak In many cases the differential equation $Lu=f$ can be written as the
limit of $L_\nu u_\nu=f$ (as $\nu\to\infty$), where {\it a
priori} estimates are available for each $L_\nu u_\nu=f$.
An {\it Ascoli-Arzela} argument is then sought for the limiting case
$Lu=f$, which needs a uniform bound for a stronger norm in order to get the
convergence of a subsequence in a weaker norm.

\medbreak As an illustration we consider the case for ${\mathbb R}$ with
$L^2$ norm as the weaker norm and $L_1^2$ (the $L^2$ norm for
derivatives up to order $1$) as the stronger norm.
The effect of a change of scale $x\to\lambda x$ is different on the two norms
$$
\displaylines{\int_{\mathbb R}\left|h\right|^2dx\,\to\,\int_{\mathbb
R}\lambda\left|h\right|^2dx,\cr \int_{\mathbb
R}\left|h^\prime\right|^2dx\,\to\,\int_{\mathbb
R}\frac{1}{\lambda}\left|h^\prime\right|^2dx.\cr}
$$
We can always make an appropriate $\nu$-dependent change of scale
$\lambda_\nu$ to make the stronger norm uniformly bounded in $\nu$.

\medbreak Scaling done in a manifold $X$ separately for ever smaller
coordinate charts is equivalent to estimating
$$\int_X |F|\left|Dh\right|^2$$ instead of
$$\int_X\left|Dh\right|^2,$$ where $F$ is a smooth
function on $X$ (and $Dh$ is the first-order differentiation of
$h$).  Here $F$ is the multiplier and describes the local rescalings of
infinitesimally small coordinate charts.

\medbreak When the first derivative $Dh$ becomes large in one direction at a point, to make $L^2_1$ norm bounded, we can enlarge the coordinate in that direction at that point. It is the same as collapsing the manifold along that direction at that point.  When we fix our sight on the manifold, $Dh$ blows up, but when we fix our sight on $Dh$, the manifold collapses.

\medbreak In the limiting situation the manifold becomes a subspace in itself. The manifold is {\it unstable}.  Before the limit is reached, it is the same manifold.  At the limit, it becomes another one.  The moduli space is not Hausdorff.  The point in the moduli space representing the manifold is not closed.  The point representing the new manifold belongs to the closure of the singleton set which represents the original manifold.

\medbreak The multiplier ideal sheaf ${\mathcal I}$ defines the subspace
into which the manifold $X$ collapses.  The structure sheaf of the subspace is
${\mathcal O}_X\left/{\mathcal I}\right.$.
From this viewpoint the multiplier ideal sheaf is known as a {\it destabilizing subsheaf}.  The subspace is the {\it destabilizing subspace}.

\medbreak If we define stability as the nonexistence of a nontrivial destabilizing subsheaf, it would just be a tautology to say that the partial differential equation is solvable if and only if we have {\it stability} of the manifold.  The challenge is to find a way to formulate this notion of stability in terms of easily verifiable conditions.  For example, there is such a good formulation for the case of solving the Hermitian-Einstein equation for the metric of a holomorphic line bundle.  The stability condition is in terms of the comparison of Chern classes of any holomorphic subbundle (or subsheaf) with the original bundle which we will explain below in (1.3).

\medbreak The one big advantage of the method of multiplier ideal sheaves is that the support of a multiplier ideal sheaf locates the set where estimates fail.  Before the advent of the theory of multiplier ideal sheaves this ``bad'' set was only investigated in the context of {\it geometric measure theory}, saying something about its Hausdorff dimension being small.  The multiplier ideal sheaf endows the ``bad set'' with geometric and analytic structures, which are inherited from the ambient manifold.  This gives us a lot of new information and new tools to work with.

\bigbreak\noindent(1.3) {\it Heuristic Discussion of Hermitian-Einstein Metrics for Stable Vector Bundles from the Viewpoint of Multiplier Ideal Sheaves as Destabilizing Subsheaves.}  Let $X$ be a compact K\"ahler manifold with K\"ahler metric $g_{i\bar j}$ and let $V$ be a holomorphic vector bundle over $X$. The problem is to determine a condition to conclude the existence of a Hermitian-Einstein metric $h_{\alpha\bar\beta}$ along the fibers of $V$ in the sense that $$\sum_{i,\bar j}g^{i\bar j}\Omega_{\alpha\bar\beta i\bar
j}=c\,h_{\alpha\bar\beta}\leqno{(1.3.1)}$$ for some constant $c$ depending on the topology of $V$, where $\Omega_{\alpha\bar\beta i\bar
j}$ is the curvature of the metric $h_{\alpha\bar\beta}$.  For our discussion we assume that $c=1$.

\medbreak The equation (1.3.1) is elliptic in the local coordinates of $X$ for the unknown $h_{\alpha\bar\beta}$.  Since the unknown
$h_{\alpha\bar\beta}$ is not a scalar unknown, we cannot conclude that we can solve the equation because of the ellipticity in the local coordinates of $X$.  When we regard $h_{\alpha\bar\beta}$ as a function on the total bundle space $V$, it becomes a scalar unknown, but the equation (1.3.1) in the local coordinates of $X$ plus the fiber coordinates of $V$ is no longer elliptic, because the equation does not involve the sum of the squares of vector fields along the fiber directions of $V$.  As described in (1.2) we can approximate the equation (1.3.1) by a sequence of equations with {\it a priori} estimates to end up with a multiplier ideal sheaf ${\mathcal I}$ (or a destabilizing subsheaf) on $V$ which is defined by the degeneracy of the sequence of solution metrics from the approximating differential equations.  The destabilizing subspace $W$ is spanned by all
nonzero eigenvectors of the limit solution $h_{\alpha\bar\beta}$.

\medbreak As integration by parts gives Kohn's algorithm in (1.1) some differential operator which when applied to multipliers produce other multipliers, in this case integration by parts yields the conclusion that $\bar\partial{\mathcal I}$ is contained in the tensor product of ${\mathcal I}$ and the bundle of $(0,1)$-forms, making $W$ a holomorphic subbundle (or a coherent subsheaf) of $V$.  As discussed in (1.2) the whole space $V$ collapses into the destabilizing subspace $W$, which inherits the holomorphic structure of $V$.

\medbreak From the equation (1.3.1) it follows that all the curvature of $V$ (after contraction by $g_{i\bar j}$) is concentrated on $W$ whose rank is strictly less than that of $V$, giving
$$
\frac{c_1(W)}{{\rm rank}\,W}>\frac{c_1(V)}{{\rm
rank}\,V}.\leqno{(1.3.2)}
$$
Here $c_1(V)$ and $c_1(W)$ mean their respective cup products with the appropriate power of the K\"ahler class of $X$.  If $V$ is assumed to be stable, the stability condition precisely stipulates the inequality direction\ \ ``\,$<$\,''\ \  in (1.3.2) for
all proper subbundles $W$ (subsheaves with torsion-free quotients)
of $V$.  It means that the limit metric $h_{\alpha\bar\beta}$ must be nondegenerate and is a Hermitian-Einstein metric for $V$.  This heuristic description explains how the condition of the stability of $V$ in terms of Chern classes guarantee the existence of a Hermitian-Einstein metric for $V$ from the viewpoint of the multiplier ideal sheaf ${\mathcal I}$ as a destabilizing subsheaf [Donaldson1985, Uhelenbeck-Yau1986, Donaldson1987, Weinkove2007].

\bigbreak\noindent(1.4) {\it Nadel's Multiplier Ideal Sheaves.}  Nadel's setting starts out with a compact complex manifold $X$ of complex dimension $n$ with the anticanonical line bundle $-K_X$ of $X$ being assumed ample.  Let $g_{i\bar j}$
be a K\"ahler metric of $X$ in the anticanonical class of $X$.  Let
$$
R_{i\bar j}=-\partial_i\partial_{\bar j}\det\left(g_{i\bar j}\right)_{1\leq i,j\leq n}
$$ be the Ricci curvature of $g_{i\bar j}$.  There is a smooth positive function $F$ on $X$ such that
$$
R_{i\bar j}-g_{i\bar j}=\partial_i\partial_{\bar j}\log F.
$$
We consider the
complex Monge-Amp\`ere equation
$$
\det\left(g_{i\bar j}+\partial_i\partial_{\bar j}\varphi\right)_{1\leq i,j\leq n}
=e^{-\varphi}F\det\left(g_{i\bar j}\right)_{1\leq i,j\leq n},\leqno{(1.4.1)}
$$
formulated by Calabi [Calabi1954a, Calabi1954b, Calabi1955] for the construction of a K\"ahler-Einstein metric of $X$.  If the equation (1.4.1) is solved, by taking $\partial\bar\partial\log$ of both sides of (1.4.1), we get
$$
-R^\prime_{i\bar j}=-\left(g^\prime_{i\bar j}-g_{i\bar j}\right)+\left(R_{i\bar j}-g_{i\bar j}\right)-R_{i\bar j}=-g^\prime_{i\bar j},
$$
(where $g^\prime_{i\bar j}=g_{i\bar j}+\partial_i\partial_{\bar j}\varphi$ and $R^\prime_{i\bar j}$ is the Ricci curvature of the K\"ahler metric $g^\prime_{i\bar j}$) and conclude that $g^\prime_{i\bar j}$ is a K\"ahler-Einstein metric of $X$.  The function $\varphi$ is a {\it K\"ahler potential perturbation} in the sense that if we locally use $\psi$ as a K\"ahler potential for $g_{i\bar j}$ so that $g_{i\bar j}=\partial_i\partial_{\bar j}\psi$, then $\varphi$ perturbs $\varphi$ to become $\psi+\varphi$ which is now a K\"ahler potential for the new metric $g^\prime_{i\bar j}$ so that $g^\prime_{i\bar j}=\partial_i\partial_{\bar j}\left(\psi+\varphi\right)$.  Continuity method is applied to solve the equation (1.4.1) by considering the solution of
$$
\det\left(g_{i\bar j}+\partial_i\partial_{\bar j}\varphi_t\right)_{1\leq i,j\leq n}
=e^{-t\varphi_t}F\det\left(g_{i\bar j}\right)_{1\leq i,j\leq n},\leqno{(1.4.2)_t}
$$
for $\varphi_t$ for $0\leq t\leq 1$, starting with $t=0$ by using [Yau1978, p.363, Theorem 1].

\medbreak The openness part of the continuity method is clear from the usual elliptic estimates and the implicit function theorem.  Nadel's multiplier ideal sheaf arises from the closedness part of the continuity method in the following way.  Suppose for some $0<t_*\leq 1$ we have a sequence $\varphi_{t_\nu}$ which satisfies $(1.4.2)_{t_\nu}$ with $t_\nu\to t_*$ monotonically strictly increasing as $\nu\to\infty$.

\medbreak \medbreak Since the first Chern class of $-K_X$, which (up to a normalizing universal constant) is represented by
$$\sum_{i,j=1}^n\left(g_{i\bar j}+\partial_i\partial_{\bar j}\varphi_t\right)\left(\frac{\sqrt{-1}}{2}dz_i\wedge d\overline{z_j}\right),\leqno{(1.4.3)_t}
$$
is independent of $t<t_*$, the $(1,1)$-form $(1.4.3)_t$ would converge weakly when $t$ goes through an appropriate sequence $t_\nu$ to $t_*$.  Let $\widehat{\varphi_t}$ be the average of $\varphi_t$ over $X$ with respect to the K\"ahler metric $g_{i\bar j}$.  Since the Green's operator for the Laplacian, with respect to the K\"ahler metric $g_{i\bar j}$, is a compact operator from the space of bounded measures on $X$ to the space of $L^1$ functions on $X$, we conclude that $\varphi_{t_\nu}-\widehat{\varphi_{t_\nu}}$ converges to some function in the $L^1$ norm for some subsequence $t_\nu$ of $t\to t_*$.

\medbreak Note that from the strict positivity of the $(1,1)$-form $(1.4.3)_t$ the Laplacian of $\varphi_t$ with respect to the K\"ahler metric $g_{i\bar j}$ is bounded from below by $-n$.  From the lower bound of Green's function we have
$$
\sup_X\varphi_t\leq\widehat{\varphi_t}+C
$$
for some constant $C$ independent of $t$ (see {\it e.g.,} [Siu1987, Chapter 3, Appendix A].

\medbreak For the other direction, by taking $-\partial\bar\partial\log$ of $(1.4.2)_t$, we get
$$
\left(R_t^\prime\right)_{i\bar j}=t\left(\left(g_t^\prime\right)_{i\bar j}-g_{i\bar j}\right)-\left(R_{i\bar j}-g_{i\bar j}\right)+R_{i\bar j}=t\left(g_t^\prime\right)_{i\bar j}+(1-t)g_{i\bar j}\geq t\left(g_t^\prime\right)_{i\bar j},
$$
where $\left(g_t^\prime\right)_{i\bar j}=g_{i\bar j}+\partial_i\partial_{\bar j}\varphi_t$ and $\left(R_t^\prime\right)_{i\bar j}$ is the Ricci curvature of the K\"ahler metric $\left(g_t^\prime\right)_{i\bar j}$.  This means that the Ricci curvature $\left(R_{t_\nu}^\prime\right)_{i\bar j}$ is bounded uniformly from below by $\left(t_*-t_\nu\right)\left(g_{t_\nu}^\prime\right)_{i\bar j}\geq\frac{t_*}{2}\left(g_{t_\nu}^\prime\right)_{i\bar j}$ for $\frac{t_*}{2}\leq t_\nu\leq t_*$. From
$$
\Delta^\prime\varphi_t=\sum_{j=1}^n\frac{\left(\varphi_t\right)_{j\bar j}}{1+\left(\varphi_t\right)_{j\bar j}}=\sum_{j=1}^n\left(1-\frac{1}{1+\left(\varphi_t\right)_{j\bar j}}\right)=n-\sum_{j=1}^n\frac{1}{1+\left(\varphi_t\right)_{j\bar j}}\leq n
$$
(evaluated with appropriate normal coordinates $z_1,\cdots,z_n$ at the point under consideration with $\Delta^\prime$ denoting the Laplacian with respect to $g^\prime_{i\bar j}$) it follows that $
\Delta^\prime\left(-\varphi_t\right)\geq-n$.  Using a Poincar\'e type inequality from lower eigenvalue estimates by a Bochner type formula and using the lower bound of the Green kernel, we get $\sup_X\left(-\varphi_{t_\nu}\right)\leq\left(n+\varepsilon\right)\sup_X\varphi_{t_\nu}+
C_\varepsilon$ for any $\varepsilon>0$ and for some constant $C_\varepsilon$ depending on $\varepsilon$ but independent of $\nu$ for $\frac{t_*}{2}\leq t_\nu\leq t_*$ (see [Siu1987, Proposition(2.2)].

\medbreak The second-order and third-order estimates used to obtain [Yau1978, p.363, Theorem 1] work also for applying the continuity method to solve $(1.4.2)_t$ for $0\leq t\leq 1$.  Alternatively the H\"older estimate for the second-order derivatives can be used instead of the third-order estimates (see {\it e.g.,} [Siu1987, Chapter 2, \S3 and \S4]).

\medbreak The obstacle in the closedness part $t\to t_*$ of the continuity method for solving $(1.4.2)_t$ occurs when $\widehat{\varphi_{t_\nu}}\to\infty$ as $\nu\to\infty$.  After multiplying $(1.4.2)_{t_\nu}$ by $e^{t_\nu\widehat{\varphi_{t_\nu}}}$ to get
$$
e^{t_\nu\widehat{\varphi_{t_\nu}}}\det\left(g_{i\bar j}+\partial_i\partial_{\bar j}\varphi_{t_\nu}\right)_{1\leq i,j\leq n}
=e^{-t_\nu\left(\varphi_{t_\nu}-\widehat{\varphi_{t_\nu}}\right)}F\det\left(g_{i\bar j}\right)_{1\leq i,j\leq n}
$$
and integrating over $X$ and taking limit as $\nu\to\infty$, we get
$$
\lim_{\nu\to\infty}\int_X e^{-t_\nu\left(\varphi_{t_\nu}-\widehat{\varphi_{t_\nu}}\right)}=\infty\leqno{(1.4.4)}
$$
when $\widehat{\varphi_{t_\nu}}\to\infty$ as $\nu\to\infty$,
because
$$
\displaylines{\int_X\det\left(g_{i\bar j}+\partial_i\partial_{\bar j}\varphi_{t_\nu}\right)_{1\leq i,j\leq n}\prod_{j=1}^n\left(\frac{\sqrt{-1}}{2}dz_j\wedge d\overline{z_j}\right)\cr=
\int_X\det\left(g_{i\bar j}\right)_{1\leq i,j\leq n}\prod_{j=1}^n\left(\frac{\sqrt{-1}}{2}dz_j\wedge d\overline{z_j}\right)=
\left(-K_X\right)^n\cr}
$$
which is independent of $t$.

\medbreak We now know that the crucial estimate in Nadel's setting is
$$
\lim_{\nu\to\infty}\int_X e^{-t_\nu\left(\varphi_{t_\nu}-\widehat{\varphi_{t_\nu}}\right)}<\infty.
$$
Since the multiplier ideal sheaf is introduced to make the crucial estimate hold after using a multiplier (in the same way as in Kohn's setting as explained in (1.1)), we introduce the multiplier ideal sheaf ${\mathcal I}$ in Nadel's setting as consisting of all holomorphic function germs $f$ on $X$ such that
$$
\lim_{\nu\to\infty}\int_U \left|f\right|^2e^{-t_\nu\left(\varphi_{t_\nu}-\widehat{\varphi_{t_\nu}}\right)}<\infty,
$$
where $U$ is an open neighborhood of the point of $X$ at which $f$ is a germ.  This multiplier ideal sheaf ${\mathcal I}$ in the sense of Nadel is defined by using a sequence of functions $\varphi_{t_\nu}-\widehat{\varphi_{t_\nu}}$ as $\nu\to\infty$ and is therefore a {\it dynamic} multiplier ideal sheaf.

\medbreak Let $\psi$ be a local plurisubharmonic function such that $g_{i\bar j}=\partial_i\partial_{\bar j}\psi$.  Since
 $$
t_\nu\left(\varphi_{t_\nu}-\widehat{\varphi_{t_\nu}}\right)+\psi=
t_\nu\left(\psi+\varphi_{t_\nu}-\widehat{\varphi_{t_\nu}}\right)+\left(1-t_\nu\right)\psi
$$
is strictly plurisubharmonic and $e^{-t_\nu\left(\varphi_{t_\nu}-\widehat{\varphi_{t_\nu}}\right)-\psi}$ is a metric for $-K_X$, it follows that ${\mathcal I}$ is a multiplier ideal sheaf for $-K_X$ and that, if (1.4.4) holds, then the multiplier ideal sheaf ${\mathcal I}$ is different from ${\mathcal O}_X$ and is therefore a nontrivial multiplier ideal sheaf for $K_X$.  The paper of Demailly-Koll\'ar [Demailly-Koll\'ar2001] uses the semi-continuity of multiplier ideal sheaves to put Nadel's multipliers in a more elegant setting.

\medbreak\noindent(1.4.5) {\it Remark on Stability Condition for Existence of K\"ahler-Einstein Metrics for Fano Manifolds.}  As discussed in (1.2), the destabilizing subspace $Y$ of $X$ defined by ${\mathcal I}$ inherits certain geometric and analytic structures from the ambient manifold $X$.  For example, we have the vanishing of $H^p\left(Y,{\mathcal O}_Y\right)$ for $p\geq 1$ by using the vanishing theorem of Nadel that $H^p\left(X,{\mathcal I}\right)=0$ for $p\geq 1$ and Kodaira's vanishing theorem $H^p\left(X,{\mathcal O}_X\right)=0$ for $p\geq 1$ and the exact long cohomology sequence of $0\to{\mathcal I}\to{\mathcal O}_X\to{\mathcal O}_Y\to 0$, because the twisting $-K_X+K_X$ of $-K_X$ by the canonical line bundle $K_X$ is the trivial line bundle.  We can regard $Y$ as the result of the collapse of $X$ as explained in (1.2) so that $Y$ inherits in some sense the Fano structure of $X$ and is itself in some sense some sort of ``Fano space''.

\medbreak When stability for a Fano manifold $X$ is defined as the impossibility of the occurrence of any nontrivial destabilizing subspace $Y$, clearly tautologically $X$ admits a K\"ahler-Einstein metric if $X$ is stable.  Unlike the case of the definition of the stability of a holomorphic vector bundle over a compact K\"ahler manifold in (1.3), there is no known easily verifiable condition which can guarantee that no such nontrivial destabilizing subspace $Y$ occurs.  From the above discussion such a condition should focus on the collapsing of $X$ into a proper subspace and not just involve the consideration of sections of vector bundles over $X$ or their subbundles or other entities defined over all of $X$.

\medbreak\noindent(1.4.6) {\it Relation Between Kohn's and Nadel's Multiplier Ideal Sheaves and their Comparison.}  Kohn's multiplier ideal sheaves and Nadel's multiplier ideal sheaves are very different in that the former consists of multipliers for the test functions (or test forms) and the latter consists of multipliers for a sequence of metrics for the anticanonical line bundle.  If we consider Nadel's vanishing theorem, Nadel's multipliers can be regarded as multipliers for the right-hand side of the $\bar\partial$ equation.  In this particular sense, since Kohn's multipliers multiply test functions (or test forms) and Nadel's multipliers multiply the right-hand side of the equation, they are in a way dual to each other.  When Kohn's multiplier ideal sheaf is nontrivial, there is no conclusion about solvability of the equation with regularity.  On the other hand, when Nadel's multiplier ideal sheaf is nontrivial in its use in Nadel's vanishing theorem, the equation can still be solved when right-hand side satisfies the condition imposed by the nontrivial multiplier ideal sheaf.

\medbreak In spite of the above fundamental differences between Kohn's and Nadel's multiplier ideal sheaves, both kinds share the following two very important features.
\begin{itemize}
\item[(i)] Both are defined by introducing multipliers into their respective crucial estimates.
\item[(ii)] Both are dynamic multiplier ideal sheaves.
\end{itemize}

\bigbreak

\bigbreak\noindent{\bf \S2. Singularity-Magnifying Complex Monge-Amp\`ere Equations.}

\bigbreak In this section we are going to discuss the construction of rational curves in Fano manifolds by producing multiplier ideal sheaves.  A trivial lemma given in (2.2) describes the kind of multiplier ideal sheaves needed for the construction of rational curves in Fano manifolds.  Such multiplier ideal sheaves cannot be constructed by using the theorem of Hirzebruch-Riemann-Roch [Hirzebruch1966] to produce the appropriate multi-valued holomorphic sections of the anticanonical line bundle $-K_X$ of the Fano manifold $X$, because of the insufficient size of the Chern number $\left(-K_X\right)^n$ for a general Fano manifold $X$ of complex dimension $n$.  This difficulty of insufficiency of the Chern number $\left(-K_X\right)^n$ is analogous to the insufficient of normal class of a curve to be deformed with two points fixed in Mori's bend-and-break argument before the introduction of Frobenius transformation in the technique of characteristic $p>0$.

\medbreak We will use appropriate complex Monge-Amp\`ere equations to produce the required multiplier ideal sheaves.  There are three kinds of complex Monge-Amp\`ere equations, which we describe as {\it singularity-reducing, singularity-neutral} and {\it singularity-magnifying}, corresponding respectively to the complex Monge-Amp\`ere equations used for the constructionn of K\"ahler-Einstein metrics for the cases of negative first Chern class, the zero first Chern class, and the positive first Chern class.

\bigbreak\noindent(2.1) {\it Three Kinds of Complex Monge-Amp\`ere Equations.}
Let $L$ be an ample line bundle on a compact complex manifold $X$ of complex dimension $n$.  Let $\sum_{i,j=1}^ng_{i\bar j}\left(\frac{\sqrt{-1}}{2}dz_i\wedge d\overline{z_j}\right)$ be a strictly positive curvature form of a smooth metric of $L$ which we are going to use as the K\"ahler form of $X$.  Let $F$ be a smooth strictly positive function on $X$.  There are the following three kinds of complex Monge-Amp\`ere equations for the unknown function $\varphi$ on the manifold $X$.
$$
\det\left(g_{i\bar j}+\partial_i\partial_{\bar j}\varphi\right)_{1\leq i,j\leq n}
=e^\varphi F\det\left(g_{i\bar j}\right)_{1\leq i,j\leq n},\leqno{(2.1.1)}
$$
$$
\det\left(g_{i\bar j}+\partial_i\partial_{\bar j}\varphi\right)_{1\leq i,j\leq n}
=F\det\left(g_{i\bar j}\right)_{1\leq i,j\leq n},\leqno{(2.1.2)}
$$
$$
\det\left(g_{i\bar j}+\partial_i\partial_{\bar j}\varphi\right)_{1\leq i,j\leq n}
=e^{-\varphi}F\det\left(g_{i\bar j}\right)_{1\leq i,j\leq n},\leqno{(2.1.3)}
$$
which are motivated respectively by complex Monge-Amp\`ere equations formulated by Calabi [Calabi1954a, Calabi1954b, Calabi1955] for K\"ahler-Einstein metrics of negative, zero, and positive first Chern class.  For equation (2.1.2) there is the following normalization condition for the function $F$
$$\int_X F\det\left(g_{i\bar j}\right)_{1\leq i,j\leq n}\prod_{j=1}^n\left(\frac{\sqrt{-1}}{2}dz_j\wedge d\overline{z_j}\right)=L^n\leqno{(2.1.4)}
$$
and the unknown function $\varphi$ is normalized by
$$
\int_X\varphi\det\left(g_{i\bar j}\right)_{1\leq i,j\leq n}\prod_{j=1}^n\left(\frac{\sqrt{-1}}{2}dz_j\wedge d\overline{z_j}\right)=0.
$$
In the original complex Monge-Amp\`ere equations formulated by Calabi for K\"ahler-Einstein metrics [Calabi1954a, Calabi1954b, Calabi1955] when $L=K_X$ in equation (2.1.1) and $L=-K_X$ in equation (2.1.3), the function $F$ is given by the condition respectively in the three cases.
$$
\displaylines{R_{i\bar j}+g_{i\bar j}=\partial_i\partial_{\bar j}\log F,\cr
R_{i\bar j}=\partial_i\partial_{\bar j}\log F,\cr
R_{i\bar j}-g_{i\bar j}=\partial_i\partial_{\bar j}\log F,\cr}
$$
with
$$
R_{i\bar j}=-\partial_i\partial_{\bar j}\det\left(g_{i\bar j}\right)_{1\leq i,j\leq n}
$$
being the Ricci curvature of $g_{i\bar j}$ so that by taking $\partial\bar\partial\log$ of both sides of each of the three complex Monge-Amp\`ere equations (2.1.1), (2.1.2), and (2.1.3) would yield respectively
$$
\displaylines{-R^\prime_{i\bar j}=g^\prime_{i\bar j}-g_{i\bar j}+\left(R_{i\bar j}+g_{i\bar j}\right)-R_{i\bar j}=g^\prime_{i\bar j},\cr
-R^\prime_{i\bar j}=R_{i\bar j}-R_{i\bar j}=0,\cr
-R^\prime_{i\bar j}=-\left(g^\prime_{i\bar j}-g_{i\bar j}\right)+\left(R_{i\bar j}-g_{i\bar j}\right)-R_{i\bar j}=-g^\prime_{i\bar j},\cr}
$$
where as above $g^\prime_{i\bar j}=g_{i\bar j}+\partial_i\partial_{\bar j}\varphi$ and $R^\prime_{i\bar j}$ is the Ricci curvature of the K\"ahler metric $g^\prime_{i\bar j}$.

\medbreak Let us first briefly discuss the different singularity behaviors of $\varphi$ in the singularity-neutral complex Monge-Amp\`ere equation (2.1.2) and the singularity-magnifying complex Monge-Amp\`ere equation (2.1.3).  We will go into these different behaviors more quantitatively in (2.3) and (2.4).  For our discussion we choose for $F$ a family $F_\varepsilon$ parametrized by $0<\varepsilon<1$ so that
$$
F_\varepsilon\bigwedge_{j=1}^n\left(\frac{\sqrt{-1}}{2}dz_j\wedge d\overline{z_j}\right)\geq
\left(\gamma\frac{\sqrt{-1}}{2}\partial\bar\partial\xi_\varepsilon\right)^n\leqno{(2.1.5)}
$$
on some coordinate chart $U$ of $X$ centered at a prescribed point $P$ of $X$,
where $\gamma$ is some positive constant and $\xi_\varepsilon$ is a smooth plurisubharmonic function on $U$ which approaches $\log|z|^2$ monotonically from above as $\varepsilon\to 0$.

\medbreak Demailly [Demailly1993] used equation (2.1.2) to produce singular metrics with strictly positive curvature current to give a partial solution to the Fujita conjecture [Fujita1987].  Because of the normalization requirement (2.1.4) for $F=F_\varepsilon$ satisfying (2.1.5), by this method the Lelong number of the singular metric $e^{-\psi-\varphi}$ of $L$ so produced (where $g_{i\bar j}=\partial_i\partial_{\bar j}\psi$ and $\varphi$ is the limit of $\varphi_\varepsilon$ for some subsequence of $\varepsilon\to 0$) is constrained to be $\leq\gamma\leq\left(L\right)^{\frac{1}{n}}$.  This is the same kind of constraint present in the method of using the theorem of Hirzebruch-Riemann-Roch [Hirzebruch1966] to obtain multi-valued holomorphic anticanonical sections to produce multiplier ideal sheaves for the construction of rational curves in Fano manifolds.  So far as the existence of rational curves in Fano manifolds by multiplier ideal sheaves is concerned, the use of equation (2.1.2) represents no advantage over the use of the theorem of Hirzebruch-Riemann-Roch [Hirzebruch1966].

\medbreak Let us now consider the use of equation (2.1.3) with $L=-K_X$ for the purpose of producing singular metrics with strictly positive curvature current to construct rational curves in Fano manifolds.  Even though we may start with $\gamma\leq\left(\left(-K_X\right)^n\right)^{\frac{1}{n}}$ and a function $F=F_\varepsilon$ satisfying (2.1.5), the factor $e^{-\varphi}$ on the right-hand side of (2.1.3) has the effect of {\it magnifying the singularity}.  We will not be able to obtain $\varphi$ as the limit of a sequence of $\varphi_\varepsilon$.  Instead, if we let $\widehat{\varphi_\varepsilon}$ be the average of $\varphi_\varepsilon$ over $X$ with respect to the K\"ahler metric $g_{i\bar j}$, then the sequence of metrics $e^{-\psi-\left(\varphi_\varepsilon-\widehat{\varphi_\varepsilon}\right)}$ of $-K_X$ produce a {\it nontrivial} multiplier ideal sheaf on $X$.  This singularity-magnifying feature of equation (2.1.3) removes for us the Chern class constraints which need to be imposed if one uses the theorem of Hirzebruch-Riemann-Roch [Hirzebruch1966] to produce rational curves in Fano manifolds.

\medbreak For the purpose of producing singular metrics with strictly positive curvature current beyond what can be accomplished by using Hirzebruch-Riemann-Roch [Hirzebruch1966], of the three complex Monge-Amp\`ere equations (2.1.1), (2.1.2), and (2.1.3), only the singularity-magnifying equation (2.1.3) is useful.  The effect of equation (2.1.2) on the singularity is neutral.  Equation (2.1.1) is even singularity-reducing.

\bigbreak\noindent(2.2) {\it Trivial Lemma.} Let $X$ be a compact complex manifold with ample anticanonical line bundle such that there exists a multiplier ideal sheaf ${\mathcal I}$ for $-K_X$ defined by a sequence of metrics $e^{-\varphi_\nu}$ for $\nu\in{\mathbb N}$ whose curvature currents have a common strictly positive lower bound.  Suppose the zero-set $Z$ of ${\mathcal I}$ is nonempty and has dimension at most one.  If ${\mathcal I}$ is not equal to the maximum ideal sheaf ${\mathfrak m}_{X,P}$ of $P$ for any point $P$ of $X$, then $Z$ has a $1$-dimensional branch whose normalization is the rational line ${\mathbb P}_1$.

\medbreak\noindent{\it Proof.}
Since the anticanonical line bundle $-K_X$ of $X$ is ample, it follows from Kodaira's vanishing theorem that $H^p\left(X,{\mathcal O}_X\right)=0$ for $p\geq 1$.  By Nadel's vanishing theorem [Nadel1990] $H^p\left(X,{\mathcal I}\right)=0$ for $p\geq 1$. From the long cohomology sequence of the short exact sequence
$$0\to{\mathcal I}\to{\mathcal O}_X\to{\mathcal O}_X\left/{\mathcal I}\right.\to 0
\leqno{(2.1.1)}
$$
it follows that $H^p\left(X,{\mathcal O}_X\left/{\mathcal I}\right.\right)=0$ for $p\geq 1$.  We now differentiate between the cases of $\dim_{\mathbb C}Z=0$ and $\dim_{\mathbb C}Z=1$.

\medbreak Suppose $\dim_{\mathbb C}Z=0$.  Since ${\mathcal I}$ is not equal to the maximum ideal sheaf ${\mathfrak m}_{X,P}$ of $P$ for any point $P$ of $X$, it follows that
$$
\dim_{\mathbb C}\Gamma\left(X,{\mathcal O}_X\left/{\mathcal I}\right.\right)\geq 2.\leqno{(2.1.2)}
$$
From $H^1\left(X,{\mathcal I}\right)=0$ the exact long cohomology sequence of (2.1.1) yields the surjectivity of
$$
\Gamma\left(X,{\mathcal O}_X\right)\to\Gamma\left(X,{\mathcal O}_X\left/{\mathcal I}\right.\right).
$$
Thus (2.1.2) implies that $\dim_{\mathbb C}\Gamma\left(X,{\mathcal O}_X\right)\geq 2$, which contradicts the fact that every holomorphic function on the compact complex manifold $X$ must be constant.  So the case $\dim_{\mathbb C}Z=0$ cannot occur and we can assume that $\dim_{\mathbb C}Z=1$.

\medbreak By choosing an appropriate $0\leq\varepsilon<1$ and an appropriately smooth metric $e^{-\psi}$ of $-K_X$ with a strictly positive curvature form and replacing
$e^{-\varphi_\nu}$ by $e^{-\left(\left(1-\varepsilon\right)\varphi_\nu+\varepsilon\psi\right)}$ for $\nu\in{\mathbb N}$, we can assume without loss of generality that $Z$ is an irreducible curve $C$ and that there exist at most a finite number of points $P_1,\cdots,P_k$ of $X$ (with possibly $k=0$) such that ${\mathcal I}$ agrees with the full ideal sheaf ${\mathcal I}_C$ of $C$ outside the points $P_1,\cdots,P_k$.

\medbreak Since the support of ${\mathcal I}_C\left/{\mathcal I}\right.$ is either empty or a finite set, it follows from the exact long cohomology sequence of the short exact sequence
$$0\to{\mathcal I}_C\left/{\mathcal I}\right.\to{\mathcal O}_X\left/{\mathcal I}\right.\to{\mathcal O}_X\left/{\mathcal I}_C\right.\to 0
$$
that $H^p\left(X,{\mathcal O}_X\left/{\mathcal I}_C\right.\right)=0$ for $p\geq 1$, which means that $H^p\left(C,{\mathcal O}_C\right)=0$ for $p\geq 1$.
Let $\pi:\tilde C\to C$ be the map for the normalization of $C$.  Then $H^p\left(\tilde C,\pi^*{\mathcal O}_C\right)=0$ for $p\geq 1$.
Since ${\mathcal O}_{\tilde C}\left/\pi^*{\mathcal O}_C\right.$
 is supported on the finite subset of $\tilde C$ which is the inverse image under $\pi$ of the singular points of $C$, it follows that $H^p\left(\tilde C,{\mathcal O}_{\tilde C}\left/\pi^*{\mathcal O}_C\right.\right)=0$ for $p\geq 1$.  From
the long cohomology sequence of the short exact sequence
$$0\to\pi^*{\mathcal O}_C\to{\mathcal O}_{\tilde C}\to{\mathcal O}_{\tilde C}\left/\pi^*{\mathcal O}_C\right.\to 0
$$
it follows that $H^p\left(X,{\mathcal O}_{\tilde C}\right)=0$ for $p\geq 1$ and $\tilde C$ is rational. Q.E.D.

\bigbreak\noindent(2.2.1) {\it Condition of Being Different from Maximum Ideal Sheaf.}  The condition that ${\mathcal I}$ in Lemma (2.2) is different from
the maximum ideal sheaf ${\mathfrak m}_{X,P}$ of $P$ for any point $P$ of $X$ can be achieved either by having $Z$ contain two points or by having $Z$ equal to a the singleton set of one point $P$ but with the subset ${\mathcal I}$ of ${\mathfrak m}_{X,P}$ strictly contained in ${\mathfrak m}_{X,P}$.

\bigbreak\noindent(2.2.2) {\it Condition of Common Strict Lower Bound for Curvature Currents.} The condition in Lemma (2.2) is essential that the curvature currents of the sequence of metrics $e^{-\varphi_\nu}$ of $-K_X$ for $\nu\in{\mathbb N}$ defining the dyanmic multiplier ideal sheaf ${\mathcal I}$ have a common strictly positive lower bound.  Suppose $X={\mathbb P}_2$ and $$s\in\Gamma\left(X,-K_X\right)=\Gamma\left({\mathbb P}_2,{\mathcal O}_{{\mathbb P}_2}(3)\right)
$$
has a nonsingular divisor with multiplicity $1$.  Then the zero-set $Z$ of the static multiplier ideal sheaf of the metric $\frac{1}{\left|s\right|^2}$ of $-K_X$ is a nonsingular elliptic curve and not the holomorphic image of a rational curve, because the curvature current of the metric $\frac{1}{\left|s\right|^2}$ does not have a strictly positive lower bound.

\medbreak Though we may not explicitly mention it for the sake of descriptional simplicity, in the rest of this note the condition of common positive lower bound for curvature currents of the sequence of metrics is assumed when we consider the nontrivial dynamic multiplier ideal sheaf produced by them for the construction of rational curves in Fano manifolds.

\bigbreak\noindent(2.3) {\it Demailly's Use of Singularity-Neutral Complex Monge-Amp\`ere Equations.}  To get results related to the Fujita conjecture, Demailly [Demailly1993] used singularity-neutral complex Monge-Amp\`ere equations to produce singular solutions for an ample line bundle $L$ over a compact complex manifold $X$ of complex dimension $n$.  Before we discuss how to use singularity-magnifying complex Monge-Amp\`ere equations to produce multiplier ideal sheaves, we first examine here Demailly's use of singularity-neutral complex Monge-Amp\`ere equations so that we can by comparison discuss more easily the singularity-magnifying effect of adding the factor $e^{-t\varphi}$ to the right-hand side of a complex Monge-Amp\`ere equation.

\medbreak Let $\sum_{i,j=1}^ng_{i\bar j}\left(\frac{\sqrt{-1}}{2}dz_i\wedge d\overline{z_j}\right)$ be a smooth strictly positive curvature form of some smooth metric of $L$, which we are going to use as the K\"ahler form of $X$ with local K\"ahler potential $\psi$ so that $g_{i\bar j}=\partial_i\partial_{\bar j}\psi$.  For $0<\varepsilon<1$ let $F_\varepsilon$ be a smooth strictly positive function on $X$.  Consider the following singularity-neutral complex Monge-Amp\`ere equation (which is obtained from the equation (2.1.2) in (2.1) by replacing $F$ in (2.1.2) by $F_\varepsilon$).
$$
\det\left(g_{i\bar j}+\partial_i\partial_{\bar j}\varphi_\varepsilon\right)_{1\leq i,j\leq n}
=F_\varepsilon\det\left(g_{i\bar j}\right)_{1\leq i,j\leq n}\leqno{(2.3.1)_\varepsilon}
$$
with $\varphi_\varepsilon$ normalized by
$$
\int_X\varphi_\varepsilon \det\left(g_{i\bar j}\right)_{1\leq i,j\leq n}\prod_{j=1}^n\left(\frac{\sqrt{-1}}{2}dz_j\wedge d\overline{z_j}\right)=0.\leqno{(2.3.2)_\varepsilon}
$$
Fix a point $P$ of $X$ and any positive number $\gamma$.  Let $U$ be a coordinate open ball neighborhood of $P$ in $X$.  We assume that $F_\varepsilon$ approaches some singular function on $X$ as $\varepsilon\to 0$.  We also assume that there exists some smooth plurisubharmonic function $\xi_\varepsilon$ on $U$ such that
\begin{itemize}\item[(i)] $\xi_\varepsilon$ approaches $\log|z|^2$ monotonically from above as $\varepsilon\to 0$,
\item[(ii)] on $U$ we have
$$
F_\varepsilon\det\left(g_{i\bar j}\right)_{1\leq i,j\leq n}\prod_{j=1}^n\left(\frac{\sqrt{-1}}{2}dz_j\wedge d\overline{z_j}\right)\geq
\left(\gamma\frac{\sqrt{-1}}{2}\partial\bar\partial\xi_\varepsilon\right)^n
$$
for $0<\varepsilon<1$, and
\item[(iii)] $F_\varepsilon$ satisfies the normalization condition
$$\int_XF_\varepsilon\det\left(g_{i\bar j}\right)_{1\leq i,j\leq n}\prod_{j=1}^n\left(\frac{\sqrt{-1}}{2}dz_j\wedge d\overline{z_j}\right)=L^n.
$$
\end{itemize}
Necessarily the constraint $\gamma^n<L^n$ occurs (or at least $\gamma^n\leq L^n$).  Demailly's use of the complex Monge-Amp\`ere equation produces singularity in the limit $\varphi$ of the solution $\varphi_\varepsilon$ as $\varepsilon\to 0$.  A conclusion of the singularity of the limit solution $\varphi$ of the solution $\varphi_\varepsilon$ as $\varepsilon\to 0$ comes from applied to a relatively compact open neighborhood $\Omega$ of $P$ in $U$ the following maximum principle of Bedford-Taylor for the complex-Monge-Amp\`ere operator given in [Bedford-Taylor1976], which is a natural generalization, from using the trace to the use of the determinant of the complex Hessian, of the usual maximum principle for the subharmonic functions based on the second derivative test of calculus.

\bigbreak\noindent(2.3.3) {\it Maximum Principle of Bedford-Taylor.}  Let $u$ and $v$ be smooth (or continuous) plurisubharmonic functions on $\bar\Omega$, where $\Omega$ is a bounded open subset of ${\mathbb C}^n$.  If
$$
u|_{\partial\Omega}\geq v|_{\partial\Omega}\ \ {\rm and\ \ }
\left(\sqrt{-1}\partial\bar\partial u\right)^n\leq \left(\sqrt{-1}\partial\bar\partial v\right)^n\ \ {\rm on\ \ }\Omega,
$$
then $u\geq v$ on $\Omega$.

\bigbreak\noindent
Since the first Chern class of $L$, which (up to a normalizing universal constant) is represented by
$$\sum_{i,j=1}^n\left(g_{i\bar j}+\partial_i\partial_{\bar j}\varphi_\varepsilon\right)\left(\frac{\sqrt{-1}}{2}dz_i\wedge d\overline{z_j}\right),
$$
is independent of $0<\varepsilon<1$, we can select a subsequence $\varepsilon_\nu\to 0$ as $\nu\to\infty$ such that the solution $\varphi_{\varepsilon_\nu}$ of $(2.3.1)_\varepsilon$ normalized by $(2.3.2)_\varepsilon$ approaches some function $\varphi$ in $L^1$ norm on $X$ as $\nu\to\infty$.  The normalization $(2.3.2)_\varepsilon$ is used only to make sure that $\varphi_{\varepsilon_\nu}$ approaches some function $\varphi$ in $L^1$ norm on $X$ as $\nu\to\infty$.  The sub-mean-value property of plurisubharmonic functions implies that there exists $C>0$ independent of $\nu$ such that
$$\psi+\varphi_{\varepsilon_\nu}\leq C+\gamma\xi_{\varepsilon_\nu}\ \ {\rm on\ }\partial\Omega\ \ {\rm for\ all}\ \nu.$$
 The
maximum principle of Bedford and Taylor (2.3.3) is now applied to
$v=\psi+\varphi_{\varepsilon_\nu}$ and $u=C+\gamma\xi_{\varepsilon_\nu}$ to yield $\psi+\varphi_{\varepsilon_\nu}\leq C+\gamma\xi_{\varepsilon_\nu}$ and  $\psi+\varphi\leq C+\gamma\log|z|$ on $\Omega$.  This means that the singularity of the limit $\varphi$ of the solution $\varphi_{\varepsilon_\nu}$ of the equation $(2.3.1)_{\varepsilon_\nu}$ normalized by $(2.3.2)_{\varepsilon_\nu}$ when $\nu\to\infty$ is no less than that of $\gamma\log|z|$ at the point $P$ with $z=0$ so far as the measurement by Lelong numbers is concerned.

\medbreak\noindent(2.3.4) {\it Singularity-Neutral Feature and Difficulty in Producing Nontrivial Multiplier Ideal Sheaves.} Though the singularity of $\varphi$ at $P$ is no less than the singularity of $\gamma\log|z|$ at the point $P$ with $z=0$, yet there is the constraint that $\gamma^n<L^n$ (or at least $\gamma^n\leq L^n$) and the order of singularity we can get for $\varphi$ is the same as what we feed into the right-hand side of the complex Monge-Amp\`ere equation.  When $L^n\leq n^n$, we have difficulty in producing a nontrivial static multiplier ideal sheaf from the metric $e^{-\varphi}$ of $L$.  When we want to use the trivial lemma (2.2) to construct rational curves in Fano manifolds $X$, the inequality $\left(-K_X\right)^n\leq n^n$ will pose the first obstacle for the use of singularity-neutral complex Monge-Amp\`ere equations.  So far as obtaining singular metrics to produce multiplier ideal sheaves is concerned, using a singularity-neutral complex Monge-Amp\`ere equation represents practically no advantage over not using it, especially when there is also some overhead cost in its use.

\bigbreak\noindent(2.4) {\it Magnification of Singularities by Singularity-Magnifying Complex Monge-Amp\`ere Equations.} We now explain how a singularity-magnifying complex Monge-Amp\`ere equation works to magnify the singularity of its solution to produce a nontrivial dynamic multiplier ideal sheaf. We use the same setting as in (2.3) simply to explain the difference in effect between the singularity-neutral and singularity-magnifying complex Monge-Amp\`ere equations.  We would like to emphasize here that for the purpose of producing rational curves on Fano manifolds some important modifications in the setting will be needed which concern some positive lower-bound condition involving the Ricci curvature and the singularity right-hand side to be formulated in (2.4.3.3) and (2.4.4.1) and explained in (2.5).  We treat the presentation of the difference between the singularity-neutral and singularity-magnifying complex Monge-Amp\`ere equations separately from the necessary modifications in order to make the singularity-magnifying argument more transparent.  The details which remain to be worked out in the analytic construction of rational curves for Fano manifolds actually lie in these modifications as will be explained in (2.5).

\medbreak We now assume that $X$ is a Fano manifold and $L=-K_X$.  We introduce a new parameter $0\leq\tau<1$ and we use the same $F_\varepsilon$ and $F$ as in (2.3), but consider the following equation $(2.4.1)_{\tau,\varepsilon}$ for the unknown $\varphi_{\tau,\varepsilon}$ instead of $(2.3.1)_\varepsilon$.
$$
\det\left(g_{i\bar j}+\partial_i\partial_{\bar j}\varphi_{\tau,\varepsilon}\right)_{1\leq i,j\leq n}
=e^{-\tau\varphi_{\tau,\varepsilon}}F_\varepsilon\det\left(g_{i\bar j}\right)_{1\leq i,j\leq n}.\leqno{(2.4.1)_{\tau,\varepsilon}}
$$
By [Yau1978, p.363, Theorem 1] for $\tau=0$ and $0<\varepsilon<1$ the equation $(2.4.1)_{\tau,\varepsilon}$ admits a solution $\varepsilon_{\tau,\varepsilon}$ and by the usual elliptic estimates and the implicit function theorem there is some $0<\tau_\varepsilon<1$ such that there is a solution $\varepsilon_{\tau,\varepsilon}$ of the equation $(2.4.1)_{\tau,\varepsilon}$ for $0\leq\tau\leq\tau_\varepsilon$ and $0<\varepsilon<1$.  We are going to prove the following simple proposition.

\medbreak\noindent(2.4.2) {\it Proposition.}  There does not exist $0<\tau_0<1$ for which there exists some monotonically decreasing sequence $\varepsilon_\nu\to 0$ as $\nu\to\infty$ such that $\tau_{\varepsilon_\nu}\geq\tau_0$ for all $\nu\in{\mathbb N}$ and the average $\widehat{\varphi_{\tau_0,\varepsilon_\nu}}$ over $X$ of $\varphi_{\tau_0,\varepsilon_\nu}$ with respect to the K\"ahler metric $g_{i\bar j}$ of $X$ is uniformly bounded for all $\nu\in{\mathbb N}$.

\medbreak\noindent{\it Proof.} Suppose the contrary and we do have such a positive number $0<\tau_0<1$ and such a monotonically decreasing sequence $\varepsilon_\nu\to 0$ with the property that
  $$
  \sup_{\nu\in{\mathbb N}}\widehat{\varphi_{\tau_0,\varepsilon_\nu}}<\infty,\leqno{(2.4.2.1)}
  $$
  where
  $$
  \widehat{\varphi_{\tau_0,\varepsilon_\nu}}=\int_X\varphi_{\tau_0,\varepsilon_\nu}\det\left(g_{i\bar j}\right)_{1\leq i,j\leq n}
  $$
  for $\nu\in{\mathbb N}$.
Since the first Chern class of $-K_X$, which (up to a normalizing universal constant) is represented by
$$\sum_{i,j=1}^n\left(g_{i\bar j}+\partial_i\partial_{\bar j}\varphi_{\tau_0,\varepsilon_\nu}\right)\left(\frac{\sqrt{-1}}{2}dz_i\wedge d\overline{z_j}\right),\leqno{(2.4.2.2)_\nu}
$$
is independent of $\nu\in{\mathbb N}$, the $(1,1)$-form $(2.4.2.2)_\nu$ would converge weakly when $\nu$ goes through an appropriate subsequence.  Since the Green's operator for the Laplacian, with respect to the K\"ahler metric $g_{i\bar j}$, is a compact operator from the space of bounded measures on $X$ to the space of $L^1$ functions on $X$, it follows from (2.4.2.1) that, by replacing the sequence $\varepsilon_\nu$ by a subsequence we can assume without loss of generality that $\varphi_{\tau_0,\varepsilon_\nu}\to\varphi$ for some $\varphi$ in $L^1$ norm on $X$.  Let $\psi$ be a local K\"ahler potential of $g_{i\bar j}$ with $g_{i\bar j}=\partial_i\partial_{\bar j}\psi$.

\medbreak The sub-mean-value property of the plurisubharmonic function $\psi+\varphi$ implies $\varphi\leq A_U$ on some neighborhood $U$ of $P$ for some constant $A_U\in{\mathbb R}$.
We now apply the maximum principle of Bedford-Taylor (2.3.3) to
$$
\left(\frac{\sqrt{-1}}{2}\partial\bar\partial\left(\psi+\varphi\right)\right)^n\geq
e^{-\tau_0 A_U}\left(\gamma\frac{\sqrt{-1}}{2}\partial\bar\partial\log\left|z\right|^2\right)^n
$$
to get
$$
\psi+\varphi\leq e^{-\frac{\tau_0 A_U}{n}}\gamma\log\left|z\right|^2+{\rm constant}\ \ {\rm on\ \ }U.\leqno{(2.4.2.3)}
$$
Since $\varphi(z)\to -\infty$ as $z\to 0$, no matter how large $B>0$ is prescribed, it follows from (2.4.2.3) that there exist some open neighborhood $W$ of $P$ in $U$ and some real number $A_W$ such that $\varphi\leq A_W\leq-B$ on $W$.  We can now conclude from the application of maximum principle of Bedford-Taylor (2.3.3) to
$$
\left(\frac{\sqrt{-1}}{2}\partial\bar\partial\left(\psi+\varphi\right)\right)^n\geq
e^{-\tau_0 A_W}\left(\gamma\frac{\sqrt{-1}}{2}\partial\bar\partial\log\left|z\right|^2\right)^n\quad{\rm on\ \ }W
$$
that
$$
\psi+\varphi\leq e^{-\frac{\tau_0 A_W}{n}}\gamma\log\left|z\right|^2+{\rm constant}\ \ {\rm on\ \ }W
$$
with $e^{-\frac{\tau_0 A_W}{n}}\geq e^{\frac{\tau_0 B}{n}}$ which goes to $\infty$ as $B\to\infty$.
This blow-up of the Lelong number at $P$ to infinity of the solution $\varphi$ gives us a contradiction, because such a Lelong number must be finite.  Q.E.D.

\bigbreak\noindent(2.4.3) {\it Nontrivial Multiplier Ideal Sheaf From Right-Hand Side of Singularity-Magnifying Complex Monge-Amp\`ere Equation Approaching Singular Limit.}  Proposition (2.4.2) says that we cannot solve the time-dependent complex Monge-Amp\`ere equation $(2.4.1)_{\tau,\varepsilon}$ for any positive time $0<\tau<1$ with uniformity in $\nu$ for any sequence of approximating $F_{\varepsilon_\nu}$ with $\varepsilon_\nu>0$ approaching $0$ as $\nu\to\infty$.   This can be interpreted in some appropriate sense as the impossibility of solving the time-dependent complex Monge-Amp\`ere equation
$$
\det\left(g_{i\bar j}+\partial_i\partial_{\bar j}\varphi_\tau\right)_{1\leq i,j\leq n}
=e^{-\tau\varphi_\tau}F\det\left(g_{i\bar j}\right)_{1\leq i,j\leq n}
$$
for the singular $F$ for any positive time $\tau>0$ no matter how small $\tau$ is.

\medbreak Because of Proposition (2.4.2) we have the following two scenarios.  The first one is that
$$
\sup_{\nu\in{\mathbb N}}\widehat{\varphi_{\tau_0,\varepsilon_\nu}}=\infty\leqno{(2.4.3.1)}
$$
even after we replace $\nu$ by a subsequence.  The second one is that, no matter how small $0<\tau_0<1$ is and how small $0<\varepsilon_0<1$ is, there exists some $0<\varepsilon_*<\varepsilon_0$ such that for some $0<\tau_*\leq\tau_0$ the closedness part of the continuity method applied to the equation $(2.4.1)_{\tau,\varepsilon_*}$ fails to produce a solution $\varphi_{\tau_*,\varepsilon_*}$ as the limit of the solution $\varphi_{\tau_\nu,\varepsilon_*}$ of $(2.4.1)_{\tau_\nu,\varepsilon_*}$ for some monotonically strictly increasing sequence $\tau_\nu\to\tau_*$ as $\nu\to\infty$.

\medbreak We now assume that we have the first scenario such that (2.4.3.1) holds.  We now multiply both sides of $(2.4.1)_{\tau_0,\varepsilon_\nu}$ by $e^{\tau_0\widehat{\varphi_{\tau_0,\varepsilon_\nu}}}$ to get $$
e^{\tau_0\widehat{\varphi_{\tau_0,\varepsilon_\nu}}}\det\left(g_{i\bar j}+\partial_i\partial_{\bar j}\varphi_{\tau_0,\varepsilon_\nu}\right)_{1\leq i,j\leq n}
=e^{-\tau_0\left(\varphi_{\tau_0,\varepsilon_\nu}-\widehat{\varphi_{\tau_0,\varepsilon_\nu}}\right)}
F_{\varepsilon_\nu}\det\left(g_{i\bar j}\right)_{1\leq i,j\leq n}.
$$
Integrating over $X$, we get
$$
\lim_{\nu\to\infty}\int_Xe^{-\tau_0\left(\varphi_{\tau_0,\varepsilon_\nu}-\widehat{\varphi_{\tau_0,\varepsilon_\nu}}\right)}
F_{\varepsilon_\nu}\det\left(g_{i\bar j}\right)_{1\leq i,j\leq n}\left(\frac{\sqrt{-1}}{2}dz_i\wedge d\overline{z_j}\right)=\infty,\leqno{(2.4.3.2)}
$$
because of (2.4.3.1) and because
$$
\int_X\det\left(g_{i\bar j}+\partial_i\partial_{\bar j}\varphi_{\tau_0,\varepsilon_\nu}\right)_{1\leq i,j\leq n}\left(\frac{\sqrt{-1}}{2}dz_i\wedge d\overline{z_j}\right)
$$
is independent of $\nu$.  We can now define the dynamic multiplier ideal sheaf ${\mathcal I}$ as consisting of all holomorphic function germs $f$ on $X$ such that
$$
\sup_\nu\int_U\left|f\right|^2e^{-\tau_0\left(\varphi_{\tau_0,\varepsilon_\nu}-\widehat{\varphi_{\tau_0,\varepsilon_\nu}}\right)}
F_{\varepsilon_\nu}\det\left(g_{i\bar j}\right)_{1\leq i,j\leq n}\left(\frac{\sqrt{-1}}{2}dz_i\wedge d\overline{z_j}\right)<\infty
$$
for some open neighborhood $U$ of the point at which $f$ is a holomorphic function germ.  We need the following condition (2.4.3.3).

\medbreak\noindent(2.4.3.3) {\it Positive Lower Bound Condition for (1,1)-Form.} For all $\nu\in{\mathbb N}$ the $(1,1)$-form
$$
-\sqrt{-1}\partial\bar\partial\log\left(e^{-\tau_0\left(\varphi_{\tau_0,\varepsilon_\nu}-\widehat{\varphi_{\tau_0,\varepsilon_\nu}}\right)}
F_{\varepsilon_\nu}\det\left(g_{i\bar j}\right)_{1\leq i,j\leq n}\right)
$$
dominates a common smooth strictly positive $(1,1)$-form on $X$ which is independent of $\nu\in{\mathbb N}$.

\medbreak\noindent
We will discuss later in (2.5) the problem of how to achieve this condition (2.4.3.3).  This condition (2.4.3.3) is needed for two reasons.  The first one is in order to obtain the coherence of the dynamic multiplier ideal sheaf ${\mathcal I}$.  The second one is in order to apply Nadel's vanishing theorem, a need as observed in (2.2.2) above.  Suppose this condition is already satisfied.  Then it follows from (2.4.3.2) that the dynamic multiplier ideal sheaf ${\mathcal I}$ for $-K_X$ is nontrivial.

\bigbreak\noindent(2.4.4) {\it Nontrivial Multiplier Ideal Sheaf From Time Approaching Critical Value in Singularity-Magnifying Factor of Complex Monge-Amp\`ere Equation.} We now continue our discussion started in (2.4.3) and consider now the second scenario so that for some $0<\tau_*\leq\tau_0$ the closedness part of the continuity method applied to the equation $(2.4.1)_{\tau,\varepsilon_*}$ fails to produce a solution $\varphi_{\tau_*,\varepsilon_*}$ as the limit of the solution $\varphi_{\tau_\nu,\varepsilon_*}$ of $(2.4.1)_{\tau_\nu,\varepsilon_*}$ for some monotonically strictly increasing sequence $\tau_\nu\to\tau_*$ as $\nu\to\infty$. We are going to use the arguments in (1.4) with appropriate adaptation to conclude that there is a nontrivial multiplier ideal sheaf for $-K_X$.  The adaptation is needed, because in (1.4) the function $F_{\varepsilon_*}$ is the difference of the $\partial\bar\partial$ potentials for the K\"ahler metric $g_{i\bar j}$ and for its Ricci curvature, whereas here $F_{\varepsilon_*}$ is a function chosen as an approximation to some singularity-mass.  We need the following condition in our adaptation of the argument of (1.4).

\medbreak\noindent(2.4.4.1) {\it Positive Lower Bound Condition for (1,1)-Form.} There exists some $\eta>0$ such that for all $\nu\in{\mathbb N}$ the following inequality for $(1,1)$-forms holds.
$$
-\sqrt{-1}\partial\bar\partial\log\left(
F_{\varepsilon_\nu}\det\left(g_{i\bar j}\right)_{1\leq i,j\leq n}\right)\geq
\eta\sum_{i,j=1}^n g_{i\bar j}\sqrt{-1}dz_i\wedge d\overline{z_j}.
$$

\medbreak\noindent We will also discuss later in (2.5) the problem of how to achieve this condition.  Now we assume that we have the condition (2.4.4.1).  Let $\left(g^\prime_{\tau,\varepsilon}\right)_{i\bar j}$ be the K\"ahler metric defined by
$\left(g^\prime_{\tau,\varepsilon}\right)_{i\bar j}=g_{i\bar j}+\partial_i\partial_{\bar j}\varphi_{\tau,\varepsilon}$ and let
$$\left(R^\prime_{\tau,\varepsilon}\right)_{k\bar\ell}=-\partial_k\partial_{\bar\ell}
\log\det\left(\left(g^\prime_{\tau,\varepsilon}\right)_{i\bar j}\right)_{1\leq i,j\leq n}$$
be its Ricci curvature.  By taking $-\partial_i\partial_{\bar j}$ of the equation
$(2.4.1)_{\tau,\varepsilon}$, we get
$$
\displaylines{\left(R^\prime_{\tau,\varepsilon}\right)_{i\bar j}=\tau\partial_i\partial_{\bar j}\varphi_{\tau,\varepsilon}-\partial_i\partial_{\bar j}\log\left(
F_{\varepsilon}\det\left(g_{k\bar\ell}\right)_{1\leq k,\ell\leq n}\right)\cr
=\tau\left(\left(g^\prime_{\tau,\varepsilon}\right)_{i\bar j}-g_{i\bar j}\right)-\partial_i\partial_{\bar j}\log\left(
F_{\varepsilon}\det\left(g_{k\bar\ell}\right)_{1\leq k,\ell\leq n}\right)\cr
\geq\tau\left(\left(g^\prime_{\tau,\varepsilon}\right)_{i\bar j}-g_{i\bar j}\right)+\eta g_{i\bar j}
=\tau\left(g^\prime_{\tau,\varepsilon}\right)_{i\bar j}+\left(\eta-\tau\right)g_{i\bar j}.\cr}
$$
When we choose $\tau_0$ at the beginning, we assume that $0<\tau_0<\eta$ so that $\eta-\tau>\eta-\tau_0$ and for $\frac{\tau_0}{2}\leq\tau\leq\tau_0$ we have the uniform lower bound
$$
\left(R^\prime_{\tau,\varepsilon}\right)_{i\bar j}\geq\frac{\tau_0}{2}
\left(g^\prime_{\tau,\varepsilon}\right)_{i\bar j}+\left(\eta-\tau_0\right)g_{i\bar j}.
$$
This is enough for the argument of (1.4) and, from the failure to produce a solution $\varphi_{\tau_*,\varepsilon_*}$ as the limit of the solution $\varphi_{\tau_\nu,\varepsilon_*}$ of $(2.4.1)_{\tau_\nu,\varepsilon_*}$ with $\tau_\nu\to\tau_*$ as $\nu\to\infty$,
we can conclude that
$$
\sup_{\nu\in{\mathbb N}}\widehat{\varphi_{\tau_\nu,\varepsilon_*}}=\infty\leqno{(2.4.4.2)}
$$
We now multiply both sides of $(2.4.1)_{\tau_\nu,\varepsilon_*}$ by $e^{\tau_\nu\widehat{\varphi_{\tau_\nu,\varepsilon_*}}}$ to get $$
e^{\tau_\nu\widehat{\varphi_{\tau_\nu,\varepsilon_*}}}\det\left(g_{i\bar j}+\partial_i\partial_{\bar j}\varphi_{\tau_\nu,\varepsilon_*}\right)_{1\leq i,j\leq n}
=e^{-\tau_\nu\left(\varphi_{\tau_\nu,\varepsilon_*}
-\widehat{\varphi_{\tau_\nu,\varepsilon_*}}\right)}
F_{\varepsilon_*}\det\left(g_{i\bar j}\right)_{1\leq i,j\leq n}.
$$
Integrating over $X$, we get
$$
\lim_{\nu\to\infty}\int_Xe^{-\tau_\nu\left(\varphi_{\tau_\nu,\varepsilon_*}
-\widehat{\varphi_{\tau_\nu,\varepsilon_*}}\right)}
F_{\varepsilon_\nu}\det\left(g_{i\bar j}\right)_{1\leq i,j\leq n}\left(\frac{\sqrt{-1}}{2}dz_i\wedge d\overline{z_j}\right)=\infty,\leqno{(2.4.4.3)}
$$
because of (2.4.4.2) and because
$$
\int_X\det\left(g_{i\bar j}+\partial_i\partial_{\bar j}\varphi_{\tau_\nu,\varepsilon_*}\right)_{1\leq i,j\leq n}\left(\frac{\sqrt{-1}}{2}dz_i\wedge d\overline{z_j}\right)
$$
is independent of $\nu$.  We can now define the dynamic multiplier ideal sheaf ${\mathcal I}$ as consisting of all holomorphic function germs $f$ on $X$ such that
$$
\sup_\nu\int_U\left|f\right|^2e^{-\tau_\nu\left(\varphi_{\tau_\nu,
\varepsilon_*}-\widehat{\varphi_{\tau_\nu,\varepsilon_*}}\right)}
F_{\varepsilon_*}\det\left(g_{i\bar j}\right)_{1\leq i,j\leq n}\left(\frac{\sqrt{-1}}{2}dz_i\wedge d\overline{z_j}\right)<\infty
$$
for some open neighborhood $U$ of the point at which $f$ is a function germ.  Because of (2.4.4.3) the dynamic multiplier ideal sheaf ${\mathcal I}$ is nontrivial.

\medbreak\noindent{(2.5)} {\it How to Handle the Two Positive Lower Bound Conditions.}  We now discuss Condition (2.4.3.3) and Condition (2.4.4.1).  Let us use a simpler set of notations by dropping some of the subscripts not directly related to our discussion and suppressing the mention of uniformity in the parameters under consideration.  We assume that we have chosen the smooth K\"ahler metric $g_{i\bar j}$ in the class of $-K_X$ such that its Ricci curvature $R_{i\bar j}$ is a strictly positive smooth closed $(1,1)$-form (for example, we can choose $R_{i\bar j}$ first and then use Yau's theorem [Yau1978, p.363, Theorem 1] to solve for $g_{i\bar j}$ in the class of $-K_X$ whose Ricci curvature is $R_{i\bar j}$).  We consider the complex Monge-Amp\`ere equation
$$
\det\left(g_{i\bar j}+\partial_i\partial_{\bar j}\varphi\right)_{1\leq i,j\leq n}
=e^{-t\varphi}F\det\left(g_{i\bar j}\right)_{1\leq i,j\leq n}.
$$
Condition (2.4.4.1) now reads that
$$
\displaylines{-\sqrt{-1}\partial\bar\partial\log\left(
F\det\left(g_{i\bar j}\right)_{1\leq i,j\leq n}\right)\cr
=-\sqrt{-1}\partial\bar\partial\log F+\sum_{i,j=1}^n R_{i\bar j}\left(\sqrt{-1}dz_i\wedge d\overline{z_j}\right)}
$$
dominates a smooth strictly positive $(1,1)$-form on $X$.
Condition (2.4.3.3) now reads that
$$
\displaylines{-\sqrt{-1}\partial\bar\partial\log\left(e^{-t\varphi}
F\det\left(g_{i\bar j}\right)_{1\leq i,j\leq n}\right)\cr
=-\sqrt{-1}\partial\bar\partial\log F+\sum_{i,j=1}^n\left(t g^\prime_{i\bar j}+\left(R_{i\bar j}-t g_{i\bar j}\right)\right)\left(\sqrt{-1}dz_i\wedge d\overline{z_j}\right)}
$$
dominates a common smooth strictly positive $(1,1)$-form on $X$.  We need only consider this for very small $t>0$.

\medbreak\noindent(2.5.1) Let $s$ be a multi-valued holomorphic section of $\delta\left(-K_X\right)$ over $X$ for some very small positive number $\delta>0$, whose divisor is a small positive rational number times a nonsingular hypersurface in $X$.  Let $h$ be any smooth metric of $-K_X$ with strictly positive curvature $\theta_{i\bar j}$ on $X$, for example,  $h=\det\left(g_{i\bar j}\right)_{1\leq i,j\leq n}$ and $\theta_{i\bar j}=R_{i\bar j}$.   We can find $t_0>0$ small enough and then $\delta_0>0$ small so that $R_{i\bar j}-tg_{i\bar j}-\delta\theta_{i\bar j}$ is strictly positive on $X$ for $0\leq t\leq t_0$ and $0\leq\delta\leq\delta_0$.  We can set $F=\frac{h^\delta}{\left|s\right|^2}$.  Then
$$
-\sqrt{-1}\partial\bar\partial\log F+\sum_{i,j=1}^n\left(t g^\prime_{i\bar j}+\left(R_{i\bar j}-t g_{i\bar j}\right)\right)\left(\sqrt{-1}dz_i\wedge d\overline{z_j}\right)
$$
dominates a common smooth strictly positive $(1,1)$-form on $X$ for $0\leq t\leq\frac{t_0}{2}$ and the two conditions are satisfied.  However, for $\delta$ small the singularity of $F$ is not enough to make $\varphi$ assume the value $-\infty$ somewhere.  Let us explain this point in more detail in order to understand the constraint which the positive lower bound conditions (2.4.3.3) and (2.4.4.1) place on the choice of the singular function $F$.

\medbreak\noindent (2.5.2) These positive lower bound conditions (2.4.3.3) and (2.4.4.1) are only needed for the singularity-magnifying complex Monge-Amp\`ere equation and are not needed for the singularity-neutral complex Monge-Amp\`ere equation.  For the case of the singularity-magnifying complex Monge-Amp\`ere equation, if we suppress the technicality of using a sequence of smooth functions to approach the singular function $F$ which is being used (2.4), the singular function $F$ times the volume form of the background K\"ahler metric $g_{i\bar j}$ is a normalizing constant times the Dirac delta at the point under consideration.  The normalizing constant is to make sure that the integral of $F$ times the volume form of the background K\"ahler metric $g_{i\bar j}$ over $X$ is equal to $\left(-K_X\right)^n$, which is the normalizing condition (2.1.4).  When $F$ (after multiplication by the volume form of $g_{i\bar j}$) is a positive constant times the Delta delta at the point under consideration, even though the positive constant is very small, the singularity of $F$ is enough to guarantee that the K\"ahler potential perturbation $\varphi$ assumes $-\infty$ at the point under consideration and the argument of singularity-magnification works.

\medbreak In the above construction of a singular $F$ satisfying the positive lower bound conditions in (2.5.1) its singularity which is only a (possibly very small) positive fractional pole-order along a divisor may not be comparable to the singularity order of a Dirac delta.  As a result the the K\"ahler potential perturbation $\varphi$ may stay locally bounded away from below from $-\infty$.  In such a case the argument of singularity-magnification cannot be applied.

\medbreak The question arises whether even for the singularity-magnifying complex Monge-Amp\`ere equation we can still use a small positive constant times the Dirac delta.  We can do it if we can represent it as the limit of functions which (after multiplication by the volume form of $g_{i\bar j}$) satisfy the positive lower-bound condition.  Can such an approximation be done, perhaps after some appropriate modifications?  The answer to such a question has not yet been explored.  An indication of a fruitful development along this line, perhaps with the analog of Dirac delta for a hypersurface instead of with the Dirac delta of a point, is the following completely trivial statement in complex dimension one, because the semipositive lower bound condition for curvature corresponds to log plurisuperharmonicity.

\medbreak\noindent (2.5.3) {\it Lemma.}  Locally the Dirac delta at the origin of ${\mathbb C}$ is the limit of log plurisuperharmonic functions.

\medbreak\noindent{\it Proof.} Let $a>0$.  Since
$$
\displaylines{
\int_{\mathbb C}\frac{a^2}{\pi\left(z\bar z+a^2\right)^2}=2a^2\int_{r=0}^\infty
\frac{rdr}{\left(r^2+a^2\right)^2}\cr=a^2\int_{s=0}^\infty
\frac{ds}{\left(s+a^2\right)^2}=a^2\left[-\frac{1}{s+a^2}\right]_{s=0}^\infty=1,\cr
}
$$
it follows from
$$
\lim_{a\to 0}\frac{a^2}{\pi\left(z\bar z+a^2\right)^2}=0\ \ {\rm for\ \ }z\not=0
$$
that the limit of $$\frac{a^2}{\pi\left(z\bar z+a^2\right)^2}$$ as $a\to 0$ is the Dirac delta at the origin.  On the other hand, $$\frac{a^2}{\pi\left(z\bar z+a^2\right)^2}$$ is log plurisuperharmonic, because of the positivity of
$$
-\partial_z\partial_{\bar z}\log\frac{a^2}{\pi\left(z\bar z+a^2\right)^2}=\frac{2a^2}{\pi\left(z\bar z+a^2\right)^2}.
$$
This computation is, of course, simply that of the curvature computation of the Fubini-Study metric of ${\mathbb P}_1$ which is Einstein.  Q.E.D.

\medbreak\noindent(2.5.4) {\it Later Starting Time for Singularity of K\"ahler Potential Perturbation.}  The attempt to use the Dirac delta is to make sure that even before we use the singularity-magnifying argument the K\"ahler potential perturbation $\varphi$ already has $-\infty$ value so that we can apply the singularity-magnifying argument to increase the Lelong number of its $\partial\bar\partial$ and get a nontrivial multiplier ideal sheaf.  The fact that even before we use the singularity-magnifying argument the K\"ahler potential perturbation $\varphi$ already has $-\infty$ value means that at time $t=0$ the K\"ahler potential perturbation $\varphi$ already has $-\infty$ value.  We can relax our requirement and demand that only at some positive time $t=t_0>0$ the K\"ahler potential perturbation $\varphi$ already has $-\infty$ value.  Instead of invoking Yau's theorem [Yau1978, p.363, Theorem 1] we solve the complex Monge-Amp\`ere equation
$$
\det\left(g_{i\bar j}+\partial_i\partial_{\bar j}\varphi\right)_{1\leq i,j\leq n}
=e^{-t_0\varphi}F\det\left(g_{i\bar j}\right)_{1\leq i,j\leq n}
$$
by the continuity method.  If we fail to reach $t=t_0$ by the continuity method we already have a nontrivial multiplier ideal sheaf for $-K_X$.  If we manage to reach $t=t_0$ in the continuity method, the process of reaching $t=t_0$ involves solving at each $t$ before $t=t_0$ the linearized form of the complex Monge-Amp\`ere equation which means inverting the (geometric nonnegative) Laplacian minus $t$ with respect to the new metric at time $t$.  The solution $\varphi$ at $t=t_0$ is obtained by continuously integrating with respect to $t$ from $t=0$ to $t=t_0$ the result obtained by inverting the operator which is equal to Laplacian minus $t$ with respect to the metric at time $t$.  We have to coordinate the choice of the background K\"ahler metric $g_{i\bar j}$ and the singular function to get the K\"ahler potential perturbation $\varphi$ to become $-\infty$ somewhere at $t=t_0$ by using appropriate {\it a priori} estimates.  The details of this technique have not yet been developed.

\medbreak\noindent(2.5.5) {\it Singularity for Right-Hand Side by Prescribing New Metric at Positive Time.}  There is another technique of handling the positive lower bound conditions (2.4.3.3) and (2.4.4.1).
It is to use a new metric at some small positive time $t=t_0>0$ to write down the right-hand side of the singularity-magnifying complex Monge-Amp\`ere equation and then get the nontrivial multiplier ideal sheaf at some later time $t>t_0$.  This technique is closely related to the discussion in (2.5.2) to consider a singular $F$ whose singularity is comparable to a Dirac delta at a point or along a divisor.

\medbreak We use the notation in (2.5.1) and use
$$
-\sqrt{-1}\partial\bar\partial\log\frac{h^{1-\delta}}{\left|s\right|^2}=-\sqrt{-1}\left(1-\delta\right)\partial
\bar\partial\log h+\sqrt{-1}\partial\bar\partial\log\left|s\right|^2
$$
as the new singular metric $g^\prime_{i\bar j}$ for $-K_X$ and then determine $F$ so that
$$
F=\frac{\ \det\left(g^\prime_{i\bar j}\right)_{1\leq i,j\leq n}\ }{\det\left(g_{i\bar j}\right)_{1\leq i,j\leq n}}
$$
with some appropriate interpretation of taking the determinant of $\left(g^\prime_{i\bar j}\right)_{1\leq i,j\leq n}$.  In this case the K\"ahler form of the singular K\"ahler metric $g^\prime_{i\bar j}$ is a closed positive $(1,1)$-current.  When we write $g^\prime_{i\bar j}=g_{i\bar j}+\partial_i\partial_{\bar j}\varphi$, the K\"ahler potential perturbation $\varphi$ has $-\infty$ values at the zero-set of $s$, but we have to worry about the appropriate definition for the Ricci tensor $R^\prime_{i\bar j}$ and its positivity.  In the first place, to appropriately define the Ricci tensor $R^\prime_{i\bar j}$, we have to approximate $g^\prime_{i\bar j}$ first by some smooth K\"ahler form $\left(g^\prime_\varepsilon\right)_{i\bar j}$ and use its Ricci tensor $\left(R^\prime_\varepsilon\right)_{i\bar j}$ in the process.  The formation of the Ricci tensor involves taking $\partial\bar\partial$ of the logarithm of the determinant of the K\"ahler metric, making the approximation process more complicated, because there is the problem of how to handle the multiplication of closed positive $(1,1)$-currents appropriately for our purpose.  There are various techniques available to smooth out closed positive $(1,1)$-currents, but here we have to worry about what happens to the Ricci curvature of the K\"ahler metrics which smooth out the closed positive $(1,1)$-current.

\medbreak We would like to inject here a remark about how to define the positivity of the Ricci curvature
without using differentiation.  It is analogous to define plurisubharmonicity by the sub-mean-value property on holomorphic disks instead of using $\partial\bar\partial$.  The Ricci curvature occurs in the Schwarz lemma of comparing volumes, and the Ricci curvatures of two K\"ahler metrics can be compared by using volume forms [Ahlfors1938, Yau1078a].  This kind of comparison is related to the maximum principle of Bedford-Taylor for the complex Monge-Amp\`ere operator in (2.3.3).  Again the details of this technique have not yet been developed.

\medbreak\noindent(2.5.6) {\it Use of K\"ahler Cone.}  Another way of producing destabilizing subsheaves is the use of boundary points of the K\"ahler cone.  Let $X$ be a Fano manifold and let ${\mathcal A}$ be the space of all smooth K\"ahler metrics on $X$ in the class $-K_X$.  From the assumptions listed below we can get a nontrivial multiplier ideal sheaf for $-K_X$ (which is a destabilizing subsheaf from the stability viewpoint of dynamic multiplier ideal sheaves introduced for crucial estimates):
\begin{itemize}
\item[(i)] some closed positive $(1,1)$-current $\omega_*$ on $X$ whose local potential assumes $-\infty$ value somewhere on $X$,

\item[(ii)] some positive number $0<\tau_0<1$, and

\item[(iii)] some sequence of elements $g^{(\nu)}=\left(g^{(\nu)}_{i\bar j}\right)_{1\leq i,j\leq n}$ in ${\mathcal A}$ for $\nu\in{\mathbb N}$
\end{itemize}
such that
\begin{itemize}
\item[(a)] the K\"ahler form of $g^{(\nu)}$ approaches $\omega_*$ weakly as $\nu\to\infty$, and

\item[(b)] for some smooth strictly positive $(1,1)$-form $\eta$ on $X$
$$\sum_{1\leq i,j\leq n}\left(R^{(\nu)}_{i\bar j}-\tau_0 g^{(\nu)}_{i\bar j}\right)\sqrt{-1}\,dz_i\wedge d\overline{z_j}\geq\eta
$$
for all $\nu\in{\mathbb N}$, where $R^{(\nu)}_{i\bar j}$ is the Ricci curvature of $g^{(\nu)}_{i\bar j}$.
\end{itemize}
This set of assumptions involves the boundary points of the K\"ahler cone of the anticanonical class.  The paper of Demailly-Paun [Demailly-Paun2004] introduced techniques concerning the boundary points of K\"ahler cones.  They used the singularity-neutral complex Monge-Amp\`ere equation.   For our purpose the corresponding techniques for the singularity-magnifying complex Monge-Amp\`ere equation have to be used instead.

\medbreak We would like to remark that Calabi obtained the three kinds of complex Monge-Amp\`ere equations (2.1.1), (2.1.2), and (2.1.3) by integrating twice the following three systems of fourth-order partial differential equations
$$
\displaylines{R^\prime_{i\bar j}+t g^\prime_{i\bar j}=\omega_{i\bar j},\cr
R^\prime_{i\bar j}=\omega_{i\bar j},\cr
R^\prime_{i\bar j}-t g^\prime_{i\bar j}=\omega_{i\bar j},\cr}
$$
with
$$
\omega_{i\bar j}=\left\{\begin{matrix}R_{i\bar j}+tg_{i\bar j}-\partial_i\partial_{\bar j}\log F\hfill\cr\cr
R_{i\bar j}-\partial_i\partial_{\bar j}\log F\hfill\cr\cr
R_{i\bar j}-tg_{i\bar j}-\partial_i\partial_{\bar j}\log F\hfill\end{matrix}\right.
$$
respectively.  That is why the singularity-magnifying complex Monge-Amp\`ere equation is to be used to study Condition (b) listed above.  Again the details of this technique have not yet been developed.

\bigbreak\noindent(2.6) {\it Dimensions of Zero-Sets of Dynamic Multiplier Ideal Sheaves.}  In order to apply the trivial lemma (2.2) to construct rational curves in Fano manifolds, we need to produce a dynamic multiplier ideal sheaf whose zero-set is of complex dimension $1$.  In the above construction of the nontrivial dynamic multiplier ideal sheaves there is no discussion about how to control their zero-sets. We are going to discuss here the question of the zero-sets of nontrivial dynamic multiplier ideal sheaves constructed by using singularity-magnifying complex Monge-Amp\`ere equations.

\medbreak In Fujita conjecture type problems, the technique of cutting down the dimension of the zero-set of a multiplier ideal sheaf is to apply the argument to the zero-set of the multiplier ideal sheaf instead of to the ambient manifold.  Here we imitate this technique.  After we have constructed a nontrivial dynamic multiplier ideal sheaves ${\mathcal I}$ constructed by using singularity-magnifying complex Monge-Amp\`ere equations, we consider the subspace $Y$ of $X$ defined by ${\mathcal I}$.  Heuristically $Y$ is obtained from $X$ by collapsing $X$ into $Y$ and, moreover, the subspace $Y$ also inherits to a certain extent the structure of $X$ and is in some sense a kind of ``Fano space''.

\medbreak We would like to apply the same method to $Y$ instead of to $X$ to produce another nontrivial dynamic multiplier ideal sheaf on $Y$ this time.  The problem is that $Y$ may be singular and it would be a great challenge to set up in a rigorous way a singularity-magnifying complex Monge-Amp\`ere equation on $Y$.  The subspace $Y$ is obtained dynamically in the sense that it can be defined in an appropriate sense as the limit of smooth closed positive $(1,1)$-forms $\alpha_\nu$ on $X$.  The Monge-Amp\`ere equation on $Y$ can be set up by using exterior product with $\alpha_\nu$.

\medbreak Dynamic multiplier ideal sheaves are defined by multipliers in a crucial
estimate.  When we use the wedge product with $\alpha_\nu$ as a substitute for the construction of a complex Monge-Amp\`ere equation on $Y$, all the estimates we consider should be uniform in $\nu$.  Details for this part are yet to be worked out.

\medbreak One more point concerning the zero-sets of dynamic multiplier ideal sheaves which we have to address is how to avoid the situation of the dynamic multiplier ideal sheaf being equal to the maximum ideal sheaf of a single point.  We need to avoid such a situation when we use the trivial lemma (2.2).  One way to address this point is to control the setup of the singularity-magnifying complex Monge-Amp\`ere equation to collapse at every step to a subspace which is only one dimension less.  The details for this have not been worked out.

\medbreak\noindent(2.7) {\it Relation with Instability of Restriction of Tangent Bundle to Destabilizing Subspace.}  We would like to make some more remarks about the two ways of producing multiplier ideal sheaves for $-K_X$ for a Fano manifold $X$ which are being compared and discussed in this note.  The two ways of producing a multiplier ideal sheaf of $-K_X$ are the following.
\begin{itemize}\item[(i)] Use holomorphic multi-valued holomorphic sections of $-K_X$ over an $n$-dimensional Fano manifold $X$ which vanish to some appropriate orders at some prescribed point $P_0$ of $X$.  Such holomorphic multi-valued holomorphic sections are obtained by using vanishing theorems and the theorem of Hirzebruch-Riemann-Roch [Hirzebruch1966], provided that the Chern number $\left(-K_X\right)^n$ is big enough (which is a very rare situation).
\item[(ii)] Use instability to produce a destabiling subsheaf which is the multiplier ideal sheaf being sought.
\end{itemize}
For the first method, though it is very rare that the Chern number $\left(-K_X\right)^n$ is big enough, yet when it is big enough the location of $P_0$ can be chosen to be any point of $X$.  On the other hand, for the second method the location of the zero-set of the multiplier ideal sheaf does not have much flexibility.  Consider the holomorphic tangent bundle $T_X$ of $X$.  In the heuristic discussion of Hermitian-Einstein metrics for stable vector bundles from the viewpoint of multiplier ideal sheaves as destabilizing subsheaves presented in (1.3), the instability of $T_X$ over $X$ means the existence over $X$ of a subbundle $W$ (or subsheaf) of $T_X$ whose normalized Chern class ({\it i.e.,} its Chern class divided by its rank) is greater than (or at least no less than) the normalized Chern class of $T_X$ (after wedging with an appropriate power of the K\"ahler class).  Here the situation is different from that of (1.3).  Here the statement about comparison of normalized Chern classes occur only for restrictions to the zero-set $Y$ of the multiplier ideal sheaf.  For simplicity of description let us assume that $Y$ is nonsingular.  Instability here means that the normalized Chern class of the tangent bundle of $Y$ on $Y$ is greater than (or at least no less than) the normalized Chern class of the restriction $T_X\big|_Y$ of $T_X$ to $Y$ (after wedging with an appropriate power of the K\"ahler class).

\medbreak When the holomorphic tangent bundle $T_X$ of $X$ is stable over $X$, its restriction $T_X\big|_C$ to a generic curve $C$ of sufficiently high degree is also stable over $C$ [Mehta-Ramanathan1982, Mehta-Ramanathan1984], but there are certain curves $C$ with the property that $T_X\big|_C$ is not stable.  An example to look at is the case of the positive-dimensional complex projective space ${\mathbb P}_n$ (with $n\geq 2$) whose tangent bundle $T_{{\mathbb P}_n}$ is stable and yet when we restrict $T_{{\mathbb P}_n}$ to a minimal rational curve $C$ we have
$$
T_{{\mathbb P}_n}\big|_C=T_C\oplus\left.{\mathcal O}_{{\mathbb P}_n}(1)\right|_C\oplus\cdots\oplus\left.{\mathcal O}_{{\mathbb P}_n}(1)\right|_C,
$$
where $T_C={\mathcal O}_{{\mathbb P}_n}(2)\big|_C$.   In this simple example, the subbundle $T_C$ destabilizes $T_{{\mathbb P}_n}\big|_C$, because the normalized Chern class of $\left.T_{{\mathbb P}_n}\right|_C$ is $\frac{n+1}{n}$ whereas the normalized Chern class of its subbundle
$T_C$ is $2>\frac{n+1}{n}$ when $n\geq 2$.  So the zero-sets of destabilizing subsheaves are in some sense rather special subvarieties whose locations do not have much flexibility.

\bigbreak\

\bigbreak\noindent{\bf APPENDIX. \ \ Key Ideas of Mori's Positive-Characteristic Proof}

\bigbreak For the benefit of analysts reading this note we highlight here the key points of Mori's argument for comparison with our approach.  We consider a one-parameter holomorphic deformation of any irreducible curve $C$ in Fano manifold $X$ of dimension $n$ with two points $P, Q$ of $C$ fixed.   If $C$ is not rational, such a deformation which ``goes around'' forces the breakup of $C$ into irreducible curves of lower genus.  If $C$ is rational, the irreducible curves obtained in this break-up are again rational curves.  This is called the {\it bend-and-break} procedure of Mori.  Starting with any irreducible curve $C$, we can continue this procedure until we get to a rational curve.  The difficulty is that such a one-parameter deformation with two points fixed may not be possible.

\medbreak For the deformation of $C$, we consider the normal bundle of its parametrizing map $f$ from the normalization $\tilde C$ of $C$ to $X$, whose determinant line bundle is $-K_X$.  The infinitesimal deformation of $f:\tilde C\to X$ with two points $P, Q$ in $\tilde C$ fixed is given by $\Gamma\left(\tilde C, {\mathcal I}_{\left\{P,Q\right\}}f^*T_X\right)$, where ${\mathcal I}_{\left\{P,Q\right\}}$ is the ideal sheaf on $\tilde C$ of the set consisting of the two points $P$ and $Q$.  The obstruction is given by $H^1\left(\tilde C, {\mathcal I}_{\left\{P,Q\right\}}f^*T_X\right)$.  In order to get a deformation, we consider the dimension of the infinitesimal deformation minus the dimension of the obstruction, which is given by the theorem of Riemann-Roch involving $\left(-K_X\right)\cdot C$ and the genus of $\tilde C$.  While a large $\left(-K_X\right)\cdot C$ is a good contribution, a large genus of $\tilde C$ is a bad contribution.  When $\tilde C$ is rational, a one-parameter deformation is always possible if $\left(-K_X\right)\cdot C\geq n+2$, because there is no bad contribution from the genus of $\tilde C$.  In general, we seek to increase the good contribution by composing $f$ with a branched cover map $\hat C\to\tilde C$ of a large number of sheets.  Unfortunately, in general because of the branching points in $\hat C\to\tilde C$, by Hurwitz's formula to compare the Euler numbers of $\hat C$ and $\tilde C$, the genus of $\hat C$ also increases, offsetting the increase in the good contribution, except when the genus of $\tilde C$ is zero, enabling us to choose $\hat C\to\tilde C$ without any branching points.

\medbreak It is at this point that positive characteristic $p>0$ plays a r\^ole by making it possible to choose a $p^m$-sheeted $\hat C\to\tilde C$ with the genus of $\hat C$ equal to that of $\tilde C$, because of the Frobenius transformation $x\to x^p$.  To illustrate the reason for this preservation of the genus, consider the case of a plane curve defined by $g(x,y)=0$ and the projection $\left(x,y\right)\mapsto x$.  The branching points on $g(x,y)=0$ for this projection are given by $\frac{\partial g}{\partial y}=0$.  In the case of characteristic $p>0$, when $y$ occurs only as $y^p$ in $g(x,y)$, the partial derivative $\frac{\partial g}{\partial y}$ is identically zero, so that every point is a branching point, which in the computation of genus, has the same effect as having no branching point.  Thus it is possible to construct a rational curve by the procedure of ``bend-and-break'' over characteristic $p>0$.

\medbreak When $X$ is in some ${\mathbb P}_N$ defined by a number of homogeneous polynomials, constructing a curve inside $X$ means adding some homogeneous polynomials to define the curve. By using modulo $p$ to go to characteristic $p>0$, we can get these additional homogeneous polynomials modulo $p$.  We can now pass to limit as $p\to\infty$ to get our desired additional homogeneous polynomials if the degrees of these polynomials modulo $p$ are bounded independent of $p$, otherwise when we pass to limit as $p\to\infty$, we may end up with infinite series of monomials instead of polynomials.  For this we simply observe that the rational curve $C$ over characteristic $p$ can be assumed to satisfy the degree bound $\left(-K_X\right)\cdot C\leq n+1$, otherwise we can break it up into rational curves of lower degree with respect to $-K_X$ by the procedure of ``bend-and-break''.

\bigbreak\noindent{\it References}

\bigbreak\noindent[Ahflors1938] Lars V. Ahlfors,
An extension of Schwarz's lemma.
{\it Trans. Amer. Math. Soc.} \textbf{43} (1938), 359--364.

\medbreak\noindent[Angehrn-Siu1995] U. Angehrn and Y.-T. Siu. Effective freeness and point separation for adjoint bundles.  {\it Invent. Math.} \textbf{122} (1995) 291--308.

\medbreak\noindent[Birkan-Cascini-Hacon-McKernan2006] C. Birkar, P.
Cascini, C. Hacon, and J. McKernan, Existence of minimal models for
varieties of log general type, arXiv:math/0610203.

\medbreak\noindent[Bedford-Taylor1976] Eric Bedford and B. A. Taylor,
The Dirichlet problem for a complex Monge-Amp\`ere equation.
{\it Invent. Math.} \textbf{37} (1976), 1--44.

\medbreak\noindent[Calabi1954a] E. Calabi, The variation of K\"ahler metrics I: The
structure of the space; II: A minimum problem, {\it Amer. Math.
Soc. Bull.} \textbf{60} (1954), Abstract  Nos. 293--294, p.168.

\medbreak\noindent[Calabi1954b] E. Calabi, The space of K\"ahler metrics, {\it
Proc. Internat. Congress Math}.  Amsterdam, 1954, Vol. 2,
pp.206--207.

\medbreak\noindent[Calabi1955] E. Calabi, On K\"ahler manifolds with vanishing
canonical class, {\it Algebraic Geometry and Topology, A Symposium
in Honor of S. Lefschetz}, Princeton Univ. Press, Princeton, 1955,
pp.78--89.

\medbreak\noindent[Carath\'eodory1909] C. Carath\'eodory, Untersuchungen
\"uber die Grundlagen der Thermodynamik. {\it Math. Ann.}  \textbf{67} (1909), 355--386

\medbreak\noindent[Chow1939]
Wei-Liang Chow,
\"Uber Systeme von linearen partiellen Differentialgleichungen erster Ordnung. {\it Math. Ann.} \textbf{117} (1939). 98--105.

\medbreak\noindent[D'Angelo1979] John P. D'Angelo, Finite type
conditions for real hypersurfaces. {\it J. Differential Geom.} \textbf{14} (1979),
59--66.
\medbreak\noindent[Demailly1993]
Jean-Pierre Demailly,
A numerical criterion for very ample line bundles.
{\it J. Diff. Geom.} \textbf{37} (1993), 323--374.

\medbreak\noindent[Demailly-Kollar2001] Jean-Pierre Demailly and J\'anos Koll\'ar, Semi-continuity of complex singularity exponents and K\"ahler-Einstein metrics on Fano orbifolds.
{\it Ann. Sci. \'Ecole Norm. Sup.} \textbf{34} (2001), 525--556.

\medbreak\noindent[Demailly-Paun2004] Jean-Pierre Demailly and Mihai Paun,
Numerical characterization of the K\"ahler cone of a compact K\"ahler manifold.
{\it Ann. of Math.} \textbf{159} (2004), 1247--1274.

\medbreak\noindent[Donaldson1985] Simon K. Donaldson,
Anti self-dual Yang-Mills connections over complex algebraic surfaces and stable vector bundles.
{\it Proc. London Math. Soc.} \textbf{50} (1985), no. 1, 1--26.

\medbreak\noindent[Donaldson1987] Simon K. Donaldson,
Infinite determinants, stable bundles and curvature.
{\it Duke Math. J.} \textbf{54} (1987), 231--247.

\medbreak\noindent[Diederich-Fornaess1978] K. Diederich and J. E. Fornaess,
Pseudoconvex domains with real-analytic boundary. {\it Ann. of
Math.} \textbf{107} (1978), 371-384.

\medbreak\noindent[Fujita1987] T.~Fujita,
On polarized manifolds whose adjoint bundles are not semipositive.
In {\it Algebraic geometry}, Sendai, 1985, pages 167--178.
North-Holland, Amsterdam, 1987.

\medbreak\noindent[Hirzebruch1966]
F. Hirzebruch,
{\it Topological methods in algebraic geometry},
3rd ed., Springer-Verlag, New York 1966

\medbreak\noindent [Kohn1979] Joseph J. Kohn, Subellipticity of the
$\bar \partial $-Neumann problem on pseudo-convex domains:
sufficient conditions. {\it Acta Math.} \textbf{142} (1979), 79--122.

\medbreak\noindent [Lu-Yau1990] Steven Shin-Yi Lu and S.-T. Yau,
Holomorphic curves in surfaces of general type.
{\it Proc. Nat. Acad. Sci. U.S.A.} \textbf{87} (1990), 80--82.

\medbreak\noindent [Mehta-Ramanathan1982] V.~B.~Mehta and A.~Ramanathan, Semistable sheaves
on projective varieties and their restriction to curves, {\it
Math. Ann.} \textbf{258} (1982),  213--224.

\medbreak\noindent [Mehta-Ramanathan1984] V.~B.~Mehta and A.~Ramanathan, Restriction of
stable sheaves and representations of the fundamental group, {\it
Invent. Math.} \textbf{77} (1984),  163--172.

\medbreak\noindent [Miyaoka1983] Yoichi Miyaoka,
Algebraic surfaces with positive indices. Classification of algebraic and analytic manifolds (Katata, 1982), 281--301,
{\it Progr. Math.} \textbf{39}, Birkhäuser Boston, Boston, MA, 1983.

\medbreak\noindent[Mori1979] S. Mori, Projective manifolds with ample tangent bundles, Ann. of Math. 111 (1979), 593-606.

\medbreak\noindent[Nadel1990]
Alan Michael Nadel,
Multiplier ideal sheaves and K\"ahler-Einstein metrics of positive scalar curvature.
{\it Ann. of Math.} \textbf{132} (1990), 549--596.

\medbreak\noindent[Schneider-Tancredi1988] Michael Schneider and Alessandro Tancredi,
Almost-positive vector bundles on projective surfaces.
{\it Math. Ann.} \textbf{280} (1988), 537--547.

\medbreak\noindent[Siu1987]
Yum-Tong Siu,
{\it Lectures on Hermitian-Einstein metrics for stable bundles and K\"ahler-Einstein metrics.}
DMV Seminar, 8. Birkh\"auser Verlag, Basel, 1987.

\medbreak\noindent[Siu2006] Y.-T. Siu, A general non-vanishing
theorem and an analytic proof of the finite generation of the
canonical ring, arXiv:math/0610740.

\medbreak\noindent[Siu2007] Y.-T. Siu, Additional explanatory notes
on the analytic proof of the finite generation of the canonical
ring, arXiv:0704.1940.

\medbreak\noindent[Siu2007] Yum-Tong Siu,
Effective Termination of Kohn's Algorithm for Subelliptic Multipliers, arXiv:0706.4113.

\medbreak\noindent[Siu2008] Yum-Tong Siu,
Finite Generation of Canonical Ring by Analytic Method (arXiv:0803.2454),
{\it J. Sci. China} \textbf{51} (2008), 481-502.

\medbreak\noindent[Siu2008a] Yum-Tong Siu, Techniques for the Analytic Proof of
the Finite Generation of the Canonical Ring (arXiv:0811.1211), to appear in the {\it Proceedings of the Conference on Current Developments in Mathematics in
Harvard University, November 16-17, 2007}.

\medbreak\noindent[Siu-Yau1980] Yum-Tong Siu and Shing-Tung Yau,
Compact K\"ahler manifolds of positive bisectional curvature.
Invent. Math. 59 (1980), 189--204.

\medbreak\noindent[Uhlenbeck-Yau1986]
Karen Uhlenbeck and Shing-Tung Yau, On the existence of Hermitian-Yang-Mills connections in stable vector bundles. {\it Comm. Pure Appl. Math.} \textbf{39} (1986), suppl., S257--S293.  {\it ibid.}  42  (1989), 703--707.

\medbreak\noindent[Weinkove2007] Ben Weinkove,
A complex Frobenius theorem, multiplier ideal sheaves and Hermitian-Einstein metrics on stable bundles.
{\it Trans. Amer. Math. Soc.} \textbf{359} (2007), no. 4, 1577--1592.

\medbreak\noindent[Yau1978]
Shing-Tung Yau,
On the Ricci curvature of a compact K\"ahler manifold and the complex Monge-Amp\`ere equation. I.
{\it Comm. Pure Appl. Math.} \textbf{31} (1978), 339--411.

\medbreak\noindent[Yau1978a] Shing Tung Yau,
A general Schwarz lemma for K\"ahler manifolds.
Amer. J. Math. 100 (1978), 197--203.

\bigbreak\noindent{\it Author's mailing address}: Department of
Mathematics, Harvard University, Cambridge, MA 02138, U.S.A.

\medbreak\noindent {\it Author's e-mail address}:
siu@math.harvard.edu

\end{document}